\documentclass[12pt]{article}

\usepackage{algorithm}
\usepackage{algorithmic}

\usepackage[english]{babel}
\usepackage{amsmath,amsthm}
\usepackage{amsfonts}
\usepackage{graphicx}
\usepackage{epsfig}
\usepackage{epstopdf}
\usepackage{url}
\usepackage{amssymb}
\usepackage{color}
\usepackage{xy} \input xy \xyoption{all}
\usepackage[%
breaklinks%
,colorlinks%
,linkcolor=red
,anchorcolor=blue%
,pagecolor=blue%
,citecolor=blue%
,bookmarks=false%
]{hyperref}

\theoremstyle{Definition}

\numberwithin{equation}{section}

\newtheorem{theorem}{Theorem}
\newtheorem{corollary}{Corollary}
\newtheorem{lemma}{Lemma}
\newtheorem{remark}{Remark}
\newtheorem{definition}{Definition}
\newtheorem{example}{Example}

\usepackage{mathrsfs}

\usepackage{algorithm}
\usepackage{algorithmic}

\def\aa{\overline{\alpha}}

\def\proj2{\mathbb{P}^2(\mathbb{K})}
\def\projtres{\mathbb{P}^3(\mathbb{K})}
\def\myV{\mathscr{V}}
\def\myW{\mathscr{W}}
\def\myS{\mathscr{S}}

\def\myB{\mathscr{B}}
\def\myC{\mathscr{C}}

\def\myG{\mathscr{G}(\projtres)}
\def\myGG{\myG^\star}

\def\cM{\mathcal{M}}
\def\cU{\mathcal{U}}
\def\K{\mathbb{K}}

\def\cT{\mathcal{T}}

\def\Res{{\mathrm{Res}}}
\def\para{\vspace{4 mm}}

\def\cV{\mathcal{V}}

\def\cM{\mathcal{M}}
\def\cU{\mathcal{U}}
\def\K{\mathbb{K}}
\def\cT{\mathcal{T}}
\def\mapdeg{\mathrm{degMap}}

\def\ot{{\overline t}}
\def\oh{{\overline h}}

\def\ox{{\overline x}}

\def\oy{{\overline y}}

\def\proj2{\mathbb{P}^2(\mathbb{K})}

\def\projtres{\mathbb{P}^3(\mathbb{K})}
\def\myV{\mathscr{V}}
\def\myW{\mathscr{W}}
\def\myS{\mathscr{S}}

\def\myF{\mathscr{F}}

\def\myB{\mathscr{B}}
\def\myC{\mathscr{C}}

\def\cE{\mathcal{E}}

\def\cU{\mathcal{U}}
\def\myG{\mathscr{G}}

\def\myGG{\myG^*}

\def\para{\vspace{2mm}}
\def\cP{{\mathcal P}}
\def\cQ{{\mathcal Q}}

\def\gcd{{\rm gcd}}
\def\lcm{{\rm lcm}}

\def\cS{{\mathcal S}}

\def\mult{\mathrm{mult}}

\def\deg{{\rm deg}}

\def\pp{{\rm PrimPart}}

\def\sol{\mathrm{SolSpace}}

\def\ot{\,{\overline t}\,}
\def\ox{\,{\overline x}\,}
\def\oy{\,{\overline y}\,}

\def\ov{\,{\overline v}\,}

\makeatletter
\newenvironment{breakablealgorithm}
  {
   \begin{center}
     \refstepcounter{algorithm}
     \hrule height.8pt depth0pt \kern2pt
     \renewcommand{\caption}[2][\relax]{
       {\raggedright\textbf{\fname@algorithm~\thealgorithm} ##2\par}%
       \ifx\relax##1\relax 
         \addcontentsline{loa}{algorithm}{\protect\numberline{\thealgorithm}##2}%
       \else 
         \addcontentsline{loa}{algorithm}{\protect\numberline{\thealgorithm}##1}%
       \fi
       \kern2pt\hrule\kern2pt
     }
  }{
     \kern2pt\hrule\relax
   \end{center}
  }
\makeatother

\begin{document}

 \title{Birational  Reparametrizations of Surfaces}

\author{Jorge Caravantes, Sonia P\'erez-D\'{\i}az and J. Rafael Sendra\\
        Universidad de Alcal\'a \\
Dpto. de F\'{\i}sica y Matem\'aticas \\
      E-28871 Madrid, Spain  \\
jorge.caravantes@uah.es, sonia.perez@uah.es, rafael.sendra@uah.es}
\date{}          
\maketitle

\begin{abstract}
Given a unirational parameterization of a surface, we present a general algorithm to determine a birational parameterization without using parameterization algorithms. Additionally, if the surface is assumed to have a birational parametrization with empty base locus, and the input parametrization is transversal, the degree of the solution is determined in advance and the dimension of the space of solutions is reduced. As a consequence, for these cases, we present a second faster algorithm.
  \end{abstract}


\section{Introduction}\label{S-intro}
The properness of parameterizations, defined by rational functions, has to do with the injective character of the rational map  they induce, between a nonempty open Zariski subset of the parameter space and the parametrized variety. It is therefore, from the point of view of applications, or of algorithmic efficiency, or even from the point of view of the theoretical analysis of algebraic varieties, a very important property.

\para

Let us be more precise. Let $\cP(t_1,\ldots,t_r)$ be a rational parameterization of a variety $\mathcal{V}$ over a field $\K$. If the induced rational map  $\cP: \K^r \rightarrow \mathcal{V}$ is not injective, the question of deciding the existence of an injective rational parametrization of $\mathcal{V}$ and, if so, of computing it, arises. The problem about existence has its answer in Luroth's Theorem for the case of dimension one (see e.g. \cite{SWP}), and in Castenuovo's rationality Theorem for the case of dimension two (see e.g. \cite{Zar1958}); in the first case, the existence of birational parametrizations (i.e. injective ones) is ensured over any field, while in the second the requirement is that the field is algebraically closed of characteristic zero.   We recall that throughout this paper we refer to generically injective rational parametrizations by either proper or birational parametrizations.

\para

Concerning to the calculation of proper parameterizations, the problem can be stated either from the implicit equations or from a parameterization. In the first case, one derives a birational parametrization using a parametrization algorithm, In the second case, the question is to compute a birational reparametrization of an input non-birational parametrization.
The case of curves (i.e. dimension one) has been analyzed by several authors and there are effective answers in both approaches (see e.g. \cite{SWP}, \cite{sonia1}). For surfaces, the problem still lacks some answers. In the implicit case, the problem is solved (see e.g. \cite{SchichoParam}). However, if the starting point is parametric, the problem is open. Although there are some partial solutions (see \cite{sonia1}, \cite{sonia}) where the problem can be reduced, in a certain sense, to the case of curves.

\para

In this paper, we focus on the case of surfaces and from the parametric point of view. We give a general algorithm (see Algorithm \ref{alg:general}), which does not require parametrization algorithms for implicitly defined surfaces, and therefore covers the open question mentioned above. To do this, in Theorem \ref{thrm:caracterizacion}, it is shown how the problem is directly related to the generic fiber of the starting parameterization. Using this fact,  solutions to the problem are determined. Additionally, we see how by imposing two additional hypotheses, both related to the base points of the parameterization, the efficiency of the general algorithm is considerably improved. This is described in Algorithm \ref{alg:generalSpecial}. These two hypotheses are: the existence of a birational parameterization without base points and the transversality of the input parameterization. The first hypothesis allows knowing in advance the degree of the solution space (see Lemma \ref{lem:grado-sin-puntos-base}) and the second allows reducing the dimension of the solution space (see Theorem \ref{theoremL(D)}).

\para

At this point in the introduction, the reader will wonder about the need to solve the problem without leaving the parametric environment; that is, the natural question that arises is why not implicitize and then use the solution provided by the parameterization algorithms. We would like to mention some reasons that, for us, justify the strategic option of looking directly for a reparameterization, if the data is given parametrically.
A first reason is that the parameterization algorithms are not particularly simple. On the other hand, developing an algorithmic solution to the problem, that does not require parameterizing techniques, supposes having a theoretical methodology to be able to address, as a future line of research, the problem for the case of rational varieties of any dimension for which a unirational parametrization is known; we recall that the problem of parameterizing in other dimensions, different from one and two, is open. In Section \ref{S-conclusions} we briefly comment on this issue.

\para

The paper is structured as follows. In Section \ref{sec-preliminaries} we introduce, through several subsections, the preliminaries of the paper on the generic fiber and the base locus of rational maps.  Section \ref{sec:theoretical-approach} is devoted to the theoretical analysis of the problem. In this section, in Theorem  \ref{thrm:caracterizacion}, we state the keys to computationally approach the properness. Section \ref{sec:computational-approach-general} deals with the development of the general algorithm. For this purpose, in Subsections \ref{subsec:fiber-P} in \ref{subsec:fiber-S}  we discuss how to effectively deal with the generic fiber of the input parametrization and of the reparametrizing functions, respectively.  In Subsection \ref{subsec-degreeS}  the degree of the reparametrizing functions are studied. Finally, in Subsection \ref{subsec:generalAlgoritm} the general algorithm is outlined. In Section \ref{sec:computational-approach-particual}, the particular case of surfaces admitting a birational parametrization with empty locus is considered. The paper ends with a section devoted to conclusions and open related problems (see Section \ref{S-conclusions}).

\para

We finish this section by introducing the main notation used throughout this paper, and stating the problem we deal with.

\para

\noindent \textbf{Notation.} Let $\K$ be an algebraically closed field of characteristic zero. $\mathbb{P}^{k}(\K)$ denotes the $k$--dimensional projective space. Furthermore, for a generically finite rational map
\[ \begin{array}{cccc}
\cM: & \mathbb{P}^{k_1}(\K) & \dashrightarrow & \mathbb{P}^{k_2}(\K) \\
     &     \oh=(h_1:\cdots:h_{k_{1}+1}) & \longmapsto & (m_1(\oh):\cdots: m_{k_2+1}(\oh)),
\end{array}
\]
where the non-zero $m_i$ are homogenous polynomials in $\oh$ of the same degree, we denote by $\deg(\cM)$ the degree $\deg_{\oh}(m_i)$, for $m_i$ non-zero, and by $\mapdeg(\cM)$ the degree of the map $\cM$; that is,  the cardinality of the generic fiber of $\cM$ (see e.g.\ \cite{Harris:algebraic}).

 Let $f\in \mathbb{L}[t_1,t_2,t_3]$ be homogeneous and non-zero, where $\mathbb{L}$ is a field extension of $\K$. Then $\myC_{\overline{\mathbb{L}}}(f)$ denotes the projective plane curve defined by $f$ over the algebraic closure of $\mathbb{L}$. When there is no risk of ambiguity, we will simply write $\myC(f)$. For $A\in \mathbb{P}^{2}(\overline{\mathbb{L}})$, we represent by
 $\mult_A(\myC(f),\myC(g))$ the multiplicity of intersection of $\myC(f)$ and $\myC(g)$ at $A$. Also, we denote by $\mult(A,\myC(f))$ the multiplicity of $\myC(f)$ at $A$.

Finally, $\myS\subset \projtres$ represents a rational projective surface and we denote by $\deg(\myS)$, the degree of $\myS$. We assume that
\begin{equation}\label{eq-P}
\cP(\ot)=\left( p_1(\ot): p_2(\ot): p_3(\ot): p_4(\ot) \right), \,\ot=(t_1,t_2,t_3),\,\, \gcd(p_1,\ldots,p_4)=1,
\end{equation}
is a fixed non-birational projective rational parametrization of $\myS$. We assume w.l.o.g. that $p_4$ is not the zero polynomial.

\para

\noindent \textbf{Birational reparametrization problem statement.} Given $\myS$ and $\cP(\ot)$ as above, determine a birational parametrization $\cQ(\ot)$ of $\myS$ as well as a rational map $\cS:\mathbb{P}^{2}(\K)\dasharrow \mathbb{P}^{2}(\K)$ such that
\begin{equation}\label{eq-sol}
 \cP(\ot)=\cQ(\cS(\ot)).
\end{equation}
In this case, we will say that $(\cQ,\cS)$ solves the birational reparametrization problem of $\cP$.

\section{Preliminaries, Notation and Problem Statement}\label{sec-preliminaries}

In this section we recall the main preliminaries to be used throughout the paper.

\subsection{The generic fiber}\label{subsec:generic-fibra}
Let $\cM$ be a generically finite rational map. That  is, $\cM$ is a rational map
$$\cM=(m_1(t_1,\ldots,t_{k_1+1}):\cdots:m_{k_2+1}(t_1,\ldots,t_{k_1+1})): \mathbb{P}^{k_1}(\K)  \dashrightarrow   \mathbb{P}^{k_2}(\K)$$
where all non-zero polynomials $m_i$ are homogenous of the same degree and
such that there exists a non-empty Zariski open subset $\Omega\subset \mathrm{im}(\cM)$ and for every $\overline{b}\in \Omega$
the cardinality of the fiber $\cM^{-1}(\overline{b})$ is invariant and finite (see e.g. \cite{Harris:algebraic} and \cite{shafa} pg. 76). We call this number the \textsf{degree of the map} and we denote it by $\mapdeg(\cM)$. We define the \textsf{generic fiber} of $\cM$ as
\[ \myF_g(\cM):=\{\overline{\alpha}\in \mathbb{P}^{k_1}(\overline{\K(\oh)}) \,|\, \cM(\overline{\alpha})=\cM(\oh) \} \]
where $\oh=(h_1,\ldots,h_{k_1+1})$ is a tuple of independent parameters and $\overline{\K(\oh)}$ is the algebraic closure of $\K(\oh)$.
 Note that $\#(\myF_g(\cM))=\mapdeg(\cM)$.

\para

Let us assume w.l.o.g. that $m_{k_2+1}$ is not the zero polynomial; if not, a change of coordinates can be applied.
Associated to $\cM$ we may consider the affine rational map
\begin{equation}\label{eq-Ma}
\cM_a=\left(\frac{m_1(t_1,\ldots,t_{k_1},1)}{m_{k_2+1}(t_1,\ldots,t_{k_1},1)},\ldots, \frac{m_{k_2}(t_1,\ldots,t_{k_1},1)}{m_{k_2+1}(t_1,\ldots,t_{k_1},1)}\right): \mathbb{A}^{k_1}(\K)  \dashrightarrow   \mathbb{A}^{k_2}(\K).
\end{equation}
Then, the generic fiber of $\cM_a$ is defined as
\[ \myF_g(\cM_a):=\{\overline{\alpha}\in \mathbb{A}^{k_1}(\overline{\K(\oh)}) \,|\, \cM_a(\overline{\alpha})=\cM_a(\oh) \} \]
where $\oh=(h_1,\ldots,h_{k_1})$ is a tuple of independent parameters and $\overline{\K(\oh)}$ is the algebraic closure of $\K(\oh)$.
 Note that $$\#(\myF_g(\cM))=\mapdeg(\cM)=\mapdeg(\cM_a)=\#(\myF_g(\cM_a)).$$
 Moreover, let us consider the dehomogenization and homogenization maps
 \begin{equation}\label{eq-DH}
 \hspace*{-7mm}\begin{array}{l}
\begin{array}{ccl}  \mathscr{D} : \{ (a_1:\cdots:a_{k_1+1}) \in \mathbb{P}^{k_1}(\K) \,|\,
 a_{k_1+1}\neq 0\} & \longrightarrow & \mathbb{A}^{k_1}(\K)\\
  (x_1:\cdots:x_{k_1+1}) & \longmapsto & \!\!\!\!\left(\dfrac{x_1}{x_{k_1+1}},\ldots,\dfrac{x_{k_1}}{x_{k_1+1}}\right),
 \end{array}
 \\
 \noalign{\para}
\begin{array}{cccl}  \mathscr{H} : & \mathbb{A}^{k_1}(\K) & \longrightarrow & \{ (a_1:\cdots:a_{k_1+1}) \in \mathbb{P}^{k_1}(\K) \,|\,
 a_{k_1+1}\neq 0\}\\
& (y_1,\cdots,y_{k_1}) & \longmapsto & (y_1:\cdots:y_{k_1},1).
  \end{array}
  \end{array}
 \end{equation}
 Taking into account that for $\overline{\alpha}\in \myF_{g}(\cM)$,  $m_{k_2+1}(\overline{\alpha})\neq 0$ because
 $m_{k_2+1}(\oh)\neq 0$, we have that, abusing notation,
 \[ \mathscr{D}(\myF_g(\cM))=\myF_g(\cM_a), \,\,  \mathscr{H}(\myF_g(\cM_a))=\myF_g(\cM).
\]
\begin{remark}\label{rem:hypothses}
In Subsections \ref{subsec:fiber-P} and \ref{subsec:fiber-S} we will deal with the question of computing or describing the generic fibre of a parametrization and/or of a dominant rational map from $\mathbb{P}^{2}(\K)$ onto $\mathbb{P}^{2}(\K)$. Some of the techniques that will be used come from   \cite{PS-grado}.
In the following we see that, for our purposes, one of the two required hypotheses in \cite{PS-grado}, can be avoided.

Given a rational affine surface parametrization
$(a_1 /a_2 ,b_1 /b_2 ,c_1 /c_2)$, in reduced form, in  Subsection 1.2. of \cite{PS-grado} two general assumptions are introduced, namely,   $\{\nabla(a_1/a_2),\nabla(b_1/b_2)\}$ must be  linearly independent as vectors in the $\K$-vector space $\K(t_1,t_2)^2$, and    $(0:1:0)$  must  belong to none of the projective curves defined by the numerators and denominators of the parametrization. Nevertheless, this second hypothesis can be omitted when dealing with the generic fiber of the parametrization, indeed: note that, taking into account Theorem 2 in \cite{PS-grado} and Proposition 1 in \cite{PDSeSi}, one may determine the degree of the rational map induced by the rational parametrization   from the polynomial $R_1$ or $R_2$, introduced in \cite{PDSeSi},  without imposing the second assumption. The underlying idea is that the resultant w.r.t. $t_2$ (similarly with $(1:0:0)$ w.r.t. $t_1$) may read wrongly the multiplicity of the point $(0:1:0)$. Nevertheless, the polynomials $R_i$ encode the coordinates of the non-constant intersection points. This second  hypothesis is  however needed for the specialization of the computation of the fibre at  a particular point (see Lemma 6 and Theorems 5 and 6 in  \cite{PS-grado}).

Similarly, the partial degrees of an implicit equation can be computed without the assumption on the point at infinity (0:1:0). More precisely, in \cite{JSC08}, Theorems 1, 2 and 4 are obtained from the results stated above and thus, if we do not specialize on the resultant, and we work with generic points, we do not need to impose that each of the  numerators and denominators of the parametrization components passes through the point at infinity (0:1:0).

In addition, the techniques in \cite{PS-grado} are easily extended to the case of dominant rational map from $\mathbb{P}^{2}(\K)$ onto $\mathbb{P}^{2}(\K)$ (see Subsection \ref{subsec:fiber-S}) by associating an auxiliary surface parametrization. As a consequence the comments above on the hypothesis on $(0:1:0)$ are also applicable.
\end{remark}

\subsection{The base locus}\label{subsec:BaseLocus}
We recall some basic notions on base points; for further information we refer to \cite{CoxPerezSendra2020}
 and   \cite{PolinomialPerezSendra2020}. The base points of a projective rational map are the points where the map is not well-defined. In our case, we will need to speak of base points of rational maps induced by surface parametrizations and/or rational maps from $\mathbb{P}^{2}(\K)$ on $\mathbb{P}^{2}(\K)$. So, we unify both cases considering a rational map
$$\cM=(m_1(\ot):\cdots:m_{k+1}(\ot)): \mathbb{P}^{2}(\K)  \dashrightarrow   \mathbb{P}^{k}(\K)$$
where all non-zero polynomials $m_i$ are homogenous of the same degree and $\gcd(m_1,\ldots,m_k)=1$. Then,   $A\in \mathbb{P}^{2}(\K)$ is called a base point of
$\cM$ if  $A\in \bigcap_{i=1}^{k+1} \myC(m_i)$. The set of all base points of $\cM$ is called the base locus of $\cM$, and we represent it by $\myB(\cM)$; note that, since $\gcd(m_1,\ldots,m_{k+1})=1$, the base locus is either empty or finite. The multiplicity of a base point of $\cM$ is the multiplicity of the point as element of  the base locus. The multiplicity of a base point can also be seen as follows. Associated to $\cM$ we introduce the polynomials
\begin{equation}\label{eq-W}
W_1(\ox,\ot):=\sum_{i=1}^{k+1} x_i\, m_i(\ot),\,\,
 W_2(\oy,\ot):=\sum_{i=1}^{k+1} y_i\, m_i(\ot),
 \end{equation}
 where $x_i, y_i$ are new variables, and we consider the corresponding projective plane curves  $\myC(W_{i})$ in $\mathbb{P}^{2}(\mathbb{F})$ where $\mathbb{F}$ is the algebraic closure of $\K(\ox,\oy)$.   Note that $\myB(\cM)\subset \myC(W_1)\cap \myC(W_2)$.
 Then, the multiplicity of a base point $A\in \myB(\cM)$ is the multiplicity of intersection of the curves $\myC(W_1)$ and $\myC(W_2)$ at $A$. We denote the multiplicity of a base point $A\in \myB(\cM)$ as
 \begin{equation}\label{eq-multBasePointP}
\mult(A,\myB(\cM)).
\end{equation}
We say $\cM$ is \textit{transversal} if, for every $A\in \myB(\cM)$, it holds that
\[ \mult(A,\myB(\cM))=\mult(A,\myC(W_1))^2. \]
 In order to check the transversality of a parametrization, one may apply Algorithm 1 presented in Section 5 in \cite{PolinomialPerezSendra2020}.

Furthermore, we introduce the notion of multiplicity of the base   locus of ${\cM}$,  denoted $\mult(\myB(\cM))$, as
\begin{equation}\label{eq-multBP}
\mult(\myB(\cM)):=\sum_{A\in \myB({\cM})} \mult(A,\myB(\cM))
\end{equation}
Note that, since   $\myB(\cP)$ is either empty of finite   the sum above well defined.

\section{Theoretical Approach}\label{sec:theoretical-approach}
In this section, we prove a characterization on the pairs $(\cQ,\cS)$ solving the birational reparametrization problem for $\cP$, where $\cP$ is as in \eqref{eq-P}. This result will be crucial for the algorithmic approach. We start with some technical lemmas.

\para

\begin{lemma}\label{lem:FS=FP}
Let $(\cQ,\cS)$ be a solution of the birational reparametrization problem. Then, $\myF_g(\cS)=\myF_g(\cP)$.
\end{lemma}
\begin{proof}
Let $\overline{\alpha}\in \myF_g(\cS)$. Then, $\cS(\overline{\alpha})=\cS(\oh)$, where $\oh=(h_1,h_2)$ is the pair of new variables. So,
$\cP(\overline{\alpha})=\cQ(\cS(\overline{\alpha}))=\cQ(\cS(\oh))=\cP(\oh)$, and hence $\overline{\alpha}\in \myF_g(\cP)$. Conversely,
let $\overline{\beta}\in \myF_g(\cP)$. Then, $\cP(\overline{\beta})=\cP(\oh)$. Thus, $\cS(\overline{\beta})=\cQ^{-1}(\cP(\overline{\beta}))=
\cQ^{-1}(\cP(\oh))=\cS(\oh)$. So, $\overline{\beta}\in \myF_g(\cS)$.
\end{proof}

\para

\begin{lemma}\label{lem:grad}
Let $\cS=(s_1:s_2:s_3):  \mathbb{P}^{2}(\K)   \dashrightarrow   \mathbb{P}^{2}(\K)$ be a generically finite rational map. Then,
$$\left\{\nabla \left( \dfrac{s_1(t_1,t_2,1)}{ s_3(t_1,t_2,1)}\right), \nabla \left(\dfrac{s_2(t_1,t_2,1)}{s_3(t_1,t_2,1)}\right)\right\}$$ are linearly independent as vectors in the $\K$-vector space $\K(\oh)^2$.
\end{lemma}
\begin{proof}
Let  $\lambda,\mu\in \K$, not both zero,   such that
$\lambda \nabla s_1+\mu \nabla s_2=\bar 0$. Integrating w.r.t. $t_1$ in $\lambda \frac{\partial s_1}{\partial t_1}+\mu \frac{\partial s_2}{\partial t_1}=0$, we get that there exists $g(t_2)$ such that $\lambda s_1+\mu s_2+g(t_2)=0$. Now, differentiating the previous equality w.r.t. $t_2$, and taking into account that $\lambda \frac{\partial s_1}{\partial t_2}+\mu \frac{\partial s_2}{\partial t_2}=0$, we get that $g(t_2)$ is a constant $k$. So, $\lambda s_1+\mu s_2=-k$. But this implies that $(s_1,s_2)$ maps $\K^2$ into the line $\lambda x+\mu y=-k$ which is impossible because the map is dominant in $\K^2$.
\end{proof}

\para

\begin{lemma}\label{lem:2a2b}
Let $\cS=(s_1:s_2:s_3):  \mathbb{P}^{2}(\K)   \dashrightarrow   \mathbb{P}^{2}(\K)$ be a generically finite rational map such that
$\myF_g(\cP)=\myF_g(\cS)$.  For $j\in \{ 1,2,3\}$, $$
\Phi_{j}(t_1,t_2):=\left(\dfrac{s_1(t_1,t_2,1)}{s_3(t_1,t_2,1)}, \dfrac{s_2(t_1,t_2,1)}{s_3(t_1,t_2,1)},\dfrac{p_j(t_1,t_2,1)}{p_4(t_1,t_2,1)}\right)$$ parametrizes an affine surface whose irreducible defining polynomial has degree one w.r.t. $z$
\end{lemma}
\begin{proof}
Let $\cS_a$ be the affine version of $\cS$,  as in \eqref{eq-Ma}.
We start proving that $\myF_g(\cS_a)=\myF_g(\Phi_j)$. Let $\overline{\alpha}\in \myF_g(\Phi_j)$. Then, $\Phi_j(\overline{\alpha})=\Phi_j(\oh)$, where $\oh=(h_1,h_2)$ is the pair of new variables. In particular $\cS_a(\overline{\alpha})=\cS_a(\oh)$. So, $\overline{\alpha}\in \myF_g(\cS_a)$. Conversely, let $\overline{\alpha}\in \myF_g(\cS_a)$. Then, $\cS_a(\overline{\alpha})=\cS_a(\oh)$. On  the other hand, since $\myF_g(\cS_a)=\myF_g(\cP_a)$, then $\cP_a(\overline{\alpha})=\cP_a(\oh)$. In particular $p_j(\overline{\alpha})/p_4(\overline{\alpha})=p_j(\oh)/p_4(\oh)$. So $\phi_j(\overline{\alpha})=\phi_{j}(\oh)$ and hence $\overline{\alpha}\in \myF_g(\Phi_j)$.  In particular, this implies that $\mapdeg(\Phi_j)=\mapdeg(\cS)$ and that $\Phi_j$ parametrizes a surface. Let  $H_j(x,y,z)$ be its irreducible defining polynomial.
Now, applying Lemma \ref{lem:grad}, we have that $\Phi_j$ satisfies the hypothesis
in \cite{JSC08}, page 120; see Remark \ref{rem:hypothses}. Applying Theorem 6 in \cite{JSC08}, we get that
\[ \deg_z(H_j)=\dfrac{\mapdeg(\cS)}{\mapdeg(\Phi_j)}=1. \]
\end{proof}

 The next theorem characterizes the solutions of the birational reparametrization problem for $\cP$.

\para

\begin{theorem}\label{thrm:caracterizacion}
Let $\cQ$ and $\cS$ be   rational maps
\[ \cQ:\mathbb{P}^{2}(\K) \dashrightarrow \myS\subset \mathbb{P}^{3}(\K), \,\,
\cS=(s_1:s_2:s_3):\mathbb{P}^{2}(\K) \dashrightarrow \mathbb{P}^{2}(\K), \]
where $s_3\neq 0$ and $\cQ(\cS)=\cP$.
The following statements are equivalent
\begin{enumerate}
\item $\cQ$ and $\cS$  solve the birational reparametrization problem for $\cP$.
\item $\myF_g(\cP)=\myF_g(\cS)$.
\end{enumerate}
Furthermore, if any of the above two equivalent conditions holds then
\begin{enumerate}
\item[(a)] For $j\in \{ 1,2,3\}$,
\begin{equation}\label{eq-Phi}
\Phi_{j}(t_1,t_2):=\left(\dfrac{s_1(t_1,t_2,1)}{s_3(t_1,t_2,1)}, \dfrac{s_2(t_1,t_2,1)}{s_3(t_1,t_2,1)},\dfrac{p_j(t_1,t_2,1)}{p_4(t_1,t_2,1)}\right)\end{equation} parametrizes an affine surface defined by an irreducible polynomial of the form
\begin{equation}\label{eq-Hj}
H_{j}(x,y,z):=A_{j,1}(x,y)-A_{j,0}(x,y) z,
\end{equation}  with $A_{j,0}$ not zero.
\item[(b)]  $\cQ(\ot)$ is the homogenization, with $t_3$ as homogenization variable, of
\begin{equation}\label{eq-T}
\cT:=\left(\dfrac{A_{1,1}(t_1,t_2)}{A_{1,0}(t_1,t_2)},\dfrac{A_{2,1}(t_1,t_2)}{A_{2,0}(t_1,t_2)},\dfrac{ A_{3,1}(t_1,t_2)}{A_{3,0}(t_1,t_2)}\right).\end{equation}
\end{enumerate}
\end{theorem}
\begin{proof} We reason with  affine coordinates; note that by hypothesis $p_4\neq 0$ and $s_3\neq 0$. Let $\cP_a$ be the affine version of $\cP$, as in \eqref{eq-Ma},  and let $\cS_a$ be the affine version of $\cS$,  as in \eqref{eq-Ma}.

 Let us see that (2) implies (1). By (2), we know that $\cS$ is generically finite, in fact $\mapdeg(\cS)=\mapdeg(\cP)$. So, $\cS_a$ is dominant in $\K^2$. Moreover, by Lemma \ref{lem:2a2b}, $\Phi_{j}$ (see \eqref{eq-Phi}) defines a surface which implicit equation has the form of $H_j$ in \eqref{eq-Hj}. Since $\cS_a$ is dominant in $\K^2$,  $A_{j,0}$ does not vanish at $\Phi_j$ and, hence, $\cT$ is well-defined (see \eqref{eq-T}).  Now, using that  $H_j(\Phi_j)=0$, we get that $\cT(\cS_a)=\cP_a$. Furthermore, since the degree is multiplicative under composition, we have that $\mapdeg(\cQ)=1$, and hence $\cQ$ is a birational affine parametrization of $\myS$. Thus, (1) holds.

 (1) implies (2) follows from  Lemma \ref{lem:FS=FP}.

 Note that the second part of the theorem, statements (a) and (b), has been shown to be consequence of (2).
\end{proof}

\para

Based on the previous result we introduce the following notion.

\para

\begin{definition}\label{def:solution-space}
Let $\cP$ be as in \eqref{eq-P}. We define the \textsf{birational reparametrization solution space}, of a fixed degree $d$, as the {either the empty set or the} set of all $d$-degree projective plane curves inducing a rational map from $\mathbb{P}^{2}(\K)$ on $\mathbb{P}^{2}(\K)$ that satisfies statement (2) in Theorem \ref{thrm:caracterizacion}. We denote it by $\sol_d(\cP)$.  Let $\sol(\cP)$ be the union of all $\sol_d(\cP)$. Note that, by Castelnuovo's Theorem, $\sol(\cP)\neq \emptyset$.

For $\overline{s}:=(s_1,s_2,s_3)\in \sol_d(\cP)^3$ such that $\gcd(s_1,s_2,s_3)=1$, we say that $ \mathbb{P}^{2}(\K)\dashrightarrow \mathbb{P}^{2}(\K); \ot \mapsto (s_1(\ot):s_2(\ot):s_3(\ot))$ is the rational map induced by $\overline{s}$.
\end{definition}

\section{The Computational Approach: the general case}\label{sec:computational-approach-general}
The computational strategy will be to determine algorithmically $d$, as well as  $\sol_d(\cP)$, where $\sol_d(\cP)\neq\emptyset$.
More precisely, let us assume that we are able to  determine $\sol_d(\cP)\neq \emptyset$. Then, for  {$\cS\in \sol_d(\cP)$}, we implicitize the surface parametrizations $\Phi_1, \Phi_2,\Phi_3$ (see Theorem \ref{thrm:caracterizacion}  \eqref{eq-Phi}) to get  $\cQ$ as in Theorem \ref{thrm:caracterizacion}  (b). Now, $(\cQ,\cS)$ is a solution of the birational reparametrization for $\cP$.
\subsection{On the generic fiber of $\cP$}\label{subsec:fiber-P}
In order to compute $\sol_d(\cP)$ we need to determine $\mapdeg(\cP)$ and moreover the generic fiber $\myF_g(\cP)$. Using Subsection \ref{subsec:generic-fibra}, we may work affinely. So, let $\cP_a$ be the affine parametrization obtained from $\cP$ as in \eqref{eq-Ma}. Let $\cP_a(t_1,t_2)$ be expressed as
\begin{equation}\label{eq-Pa}
\cP_a(t_1,t_2)=\left(\dfrac{P_1(t_1,t_2)}{Q_1(t_1,t_2)},\dfrac{P_2(t_1,t_2)}{Q_2(t_1,t_2)},\dfrac{P_3(t_1,t_2)}{Q_3(t_1,t_2)} \right)
\end{equation}
where the rational functions are in reduced form. We show how to describe the points of $\myF_g(\cP_a)$. We consider the polynomials
\begin{equation}\label{eq-G}
 G_ i=P_i(h_1,h_2) Q_i(t_1,t_2)-P_i(t_1,t_2)Q_i(h_1,h_2), \,\,i\in \{1,2,3\}.
 \end{equation}
 Furthermore, let $W:=w \cdot \lcm(Q_1,Q_2,Q_3)-1$ where $w$ is a new variable. Also, we consider the projection
 $$
 \begin{array}{ccc}\pi: \overline{\K(h_1,h_2)}^3 & \rightarrow &  \overline{\K(h_1,h_2)}^2 \\
  (t_1,t_2,w) &\mapsto & (t_1,t_2)
  \end{array}.
 $$

Then, we have the following lemma.
\para

\begin{lemma}\label{lemma:fibra1}
$\myF_g(\cP_a) = \overline{\pi(\mathbb{V}_{\overline{\K(h_1,h_2)}}(G_1,G_2,G_3,W))}$.
\end{lemma}
\begin{proof} Let $A:=\lcm(Q_1,Q_2,Q_3)$ and $\mathscr{W}:=\mathbb{V}_{\overline{\K(h_1,h_2)}}(G_1,G_2,G_3,W)$.

Let   $\overline{\alpha}\in \myF_g(\cP_a)\subset \overline{\K(h_1,h_2)}^2$. Then $\cP_a(\overline{\alpha})=\cP_a(h_1,h_2)$. Furthermore, $\cP_a(\overline{\alpha})$ is well defined and hence $A(\aa)\neq 0$.   Therefore, $(\aa,1/A(\aa))\in \mathscr{W}$. So, $\pi((\aa,1/A(\aa)))=\aa\in \pi(\mathscr{W})\subset \overline{\pi(\mathscr{W})}$.

Conversely, let $(a,b,c)\in \mathscr{W}$. Then   $G_i(h_1,h_2,a,b)=0$ for all $i\in \{1,2,3\}$, and $c \,A(a,b)=1$. So, $A(a,b)\neq 0$. Therefore, $\cP_a(a,b)$ is well defined and $\cP_a(a,b)=\cP(h_1,h_2)$. Thus, $(a,b)=\pi(a,b,c)\in \myF_g(\cP_a)$. \end{proof}

   \para

 Therefore, using the elimination property of Gr\"obner bases and the  Closure Theorem (see e.g. \cite{Cox1} page 125, and  \cite{Winkler} page 192),  if $\tilde{\mathscr{K}}$ is a Gr\"obner basis, w.r.t. lexicographic ordering with $t_1<t_2<w$, of the ideal
 $\tilde{\mathrm{J}}:=<G_1,G_2,G_3,W> \subset \K(h_1,h_2)[t_1,t_2,w]$, and let  $\mathrm{J}:=\tilde{\mathscr{K}} \cap \K(h_1,h_2)[t_1,t_2]$. Then
 \begin{equation}\label{eq-fibraGB}
\mathbb{V}_{\overline{\K(h_1,h_2)}}(\mathrm{J})=\myF_g(\cP_a).
 \end{equation}
 Furthermore, if   $\tilde{\mathscr{K}}$   is minimal then
$\tilde{\mathscr{K}}=\{k_{1,1}, k_{2,1},\ldots,k_{2,k_2}, k_{3,1},\ldots,k_{3,k_3}\}$
 where  $k_{1,1}\in \K(h_1,h_2)[t_1]$,
 $k_{2,j}\in \K(h_1,h_2)[t_1,t_2]$ and $k_{3,j}\in \K(h_1,h_2)[t_1,t_2,w]$  (see e.g. \cite{Winkler} page 194). Thus, $\myF_g(\cP_a)$ is described by
 \begin{equation}\label{eq-Generadores}
  \mathscr{K}=\{k_{1,1}, k_{2,1},\ldots,k_{2,k_2}\}.
 \end{equation}
 That is
 \begin{equation}\label{eq-F-BG}
 \myF_g(\cP_a)=\mathbb{V}(\mathscr{K})=\{ \aa\in \overline{\K( {h_1,h_2})}^2 \, | \, k(\aa)=0 \,\,\text{for $k\in \mathscr{K}$}\}.
 \end{equation}

 On the other hand, since  $\myF_g(\cP_a)$ is zero-dimensional, we may assume, maybe after a linear change  of $\{t_1,t_2\}$, that  $\mathrm{J}$ is in general position w.r.t. $t_1$ (see   \cite{Winkler} page  194). So, if we work with $\sqrt{\mathrm{J}}$, we may apply the Shape Lemma (see \cite{Winkler} page  195). More precisely, the normal reduced Gr\"obner basis of $\sqrt{\mathrm{J}}$, w.r.t. the lexicographic order with $t_1<t_2$,
 is of the form
\begin{equation}\label{eq-GBbuena}
\mathscr{G}:= \{ u(t_1),t_2-v(t_1)  \},
\end{equation}
 with $u$ square-free and $\deg(v)<\deg(u)$.  Therefore
 \begin{equation}\label{eq-F-BGBuena}
\myF_g(\cP_a)=\mathbb{V}(\mathscr{G})=\{ \aa:=(\alpha_1,v(\alpha_1))\in \overline{\K( {h_1,h_2})}^2 \, | \, u(\alpha_1)=0 \}.
 \end{equation}

 In order     to determine $\sqrt{\mathrm{J}}$ we can, for instance, use Seidenberg Lemma (see e.g. \cite{seidenberg} or \cite{laplagne}). More precisely, if $\mathrm{J}\cap \K(h_1,h_2) [t_1]=<f(t_1)>$ and $\mathrm{J}\cap  \K(h_1,h_2) [t_2]=<g(t_2)>$, then  (see \eqref{eq-F-BG})
\begin{equation}\label{eq-rad} \sqrt{\mathrm{J}}=<k_{1,1}, k_{2,1},\ldots,k_{2,k_2} , \tilde{f}, \tilde{g}>,
\end{equation}
where $\tilde{f}= f/\gcd(f,f')$ and  $\tilde{g}= g/\gcd(g,g')$.


In \cite{PS-grado}, the authors introduce two polynomials whose roots describe, respectively, the coordinates of the elements in $\myF_g(\cP_a)$, and hence its degree provides $\mapdeg(\cP)$. More precisely, let us assume that  $\{\nabla P_1/Q_1,\nabla  P_2/Q_2 \}$ are linearly independent as vectors in the $\K$ vector space $\K(t_1,t_2)^2$, see Remark \ref{rem:hypothses}. Observe that this condition can be achieved after a $(x,y,z)$-affine change of coordinates. In addition   $\mapdeg$ will not change and $\myF_g(\cP_a)$ will be obtained applying the corresponding inverse of the linear transformation. So, we  assume that $\cP_a$ does indeed satisfy the hypothesis. In this situation, let
\begin{equation}\label{eq-R1R2old}
\left\{\begin{array}{l} \hat{R}_1(h_1,h_2,t_1)=\pp_{\{h_1,h_2\}}(\mathrm{Content}_{Z}(\Res_{t_2}(G_1,G_2+Z G_3))) \\
\noalign{\medskip}
\hat{R}_2(h_1,h_2,t_2)=\pp_{\{h_1,h_2\}}(\mathrm{Content}_{Z}(\Res_{t_1}(G_1,G_2+Z G_3))).\end{array} \right.
\end{equation}
We observe that computing $\pp_{\{h_1,h_2\}}$ in \eqref{eq-R1R2old} we are avoiding the base points  of $\cP(\ot)$. So, the polynomials above can be simplified by also avoiding the base points of $\cP(\oh)$. That is, we introduce the polynomials
\begin{equation}\label{eq-R1R2}
\left\{\begin{array}{l}  {R}_1(h_1,h_2,t_1)=\pp_{t_1}(\hat{R}_1) \\
\noalign{\medskip}
{R}_2(h_1,h_2,t_2)=\pp_{t_2}(\hat{R_2}).\end{array} \right.
\end{equation}
Then, it holds that
\begin{equation}\label{eq-degree}
 \mapdeg(\cP)=\deg_{t_1}({R}_1)=\deg_{t_2}({R}_2).
 \end{equation}
Moreover, if ${R}_1, {R}_2$ are considered as polynomials in $\K(h_1,h_2)[t_1], \K(h_1,h_2)[t_2]$, respectively, the roots of ${R}_1$ (resp. of ${R}_2)$ are the first coordinates (resp. second coordinates) of the elements in the fiber.  In addition, note that $f$ and $g$ involved in  \eqref{eq-rad} are $R_1,$ and $R_2$, respectively.

\subsection{On the generic fiber of $\cS$}\label{subsec:fiber-S}
The reasonings in Subsection \ref{subsec:fiber-P} can be analogously performed to compute the generic fiber of rational maps from $\proj2$ onto $\proj2$.  
In the sequel, we see how this can be reduced to the case in Subsection \ref{subsec:fiber-P}. More precisely, let $\cS=(s_1:s_2:s_3): \proj2 \rightarrow \proj2$ be  a generically finite rational map,  where we assume w.o.l.g. that $s_3$ is not zero.  We associate to $\cS$ the map
\[ \begin{array}{lccc}
\tilde{\cS}: & \proj2 & \dashrightarrow & \mathbb{P}^{3}(\K) \\
& \ot &\longmapsto & (s_1(\ot):s_2(\ot):s_3(\ot): s_3(\ot)).
\end{array}
\]
Observe that $\tilde{\cS}(\proj2)$ is dense in $\{(x_1:\cdots :x_4)\in \mathbb{P}^3(\K) \,|\, x_3=x_4\}$. So, $\tilde{\cS}$ parametrizes the plane $x_3=x_4$.
Let  (see \eqref{eq-Ma})
\[
\tilde{\cS}_a=\left(\dfrac{A_1(t_1,t_2)}{B_1(t_1,t_2)},\dfrac{A_2(t_1,t_2)}{B_2(t_1,t_2)}, 1\right) ,
\]
where the fractions are in reduced form.
Let us check that $\tilde{\cS}_a$ satisfies {the  hypothesis in Subsection 1.2. in \cite{PS-grado}; see Remark \ref{rem:hypothses}.} We need to ensure the hypothesis requiring that $\{\nabla (A_1/B_1),\nabla (A_2/B_2)\}$ are linearly independent as vectors in the $\K$ vector space $\K(t_1,t_2)^2$. However, if they were linearly dependent, reasoning as in Subsection \ref{subsec:fiber-P}, we would get that the image of $\cS$ would be included in a line which is a contradiction with the fact that $\cS(\proj2)$ is dense in $\proj2$.  In this situation, we consider the polynomials
 \begin{equation}\label{eq-GS}
 \tilde{G}_i := A_i(h_1,h_2)B_i(t_1,t_2)-A_i(t_1,t_2)B_i(h_1,h_2),\,\,\,i\in \{1,2\}.
 \end{equation}
and
\begin{equation}\label{eq-tildeRS}
  \left\{\begin{array}{l}
  \tilde{R}_1=\pp_{t_1}(\pp_{\{h_1,h_2\}}(\Res_{t_2}(\tilde{G}_1,\tilde{G}_2))),\\  \noalign{\medskip}
  \tilde{R}_2=\pp_{t_2}(\pp_{\{h_1,h_2\}}(\Res_{t_1}(\tilde{G}_1,\tilde{G}_2))).
  \end{array}\right.
 \end{equation}
 Then
 \[ \mapdeg(\cS)=\deg_{t_1}(\tilde{R}_1)=\deg_{t_2}(\tilde{R}_2). \]

\subsection{On the degree of $\cS$}\label{subsec-degreeS}
The first problem that we find, in order to determine the solution space, is to know  $d:=\deg(\cS)$ such that $\sol_{d}(\cP)\neq \emptyset$.

\para

\begin{lemma}\label{lem:grado-general}
Let $\cS\in \sol(\cP)$ then $$\deg(\cS)=\sqrt{\mapdeg(\cP)+\mult(\myB(\cS))}\geq \lceil \sqrt{\mapdeg(\cP)} \rceil.$$
\end{lemma}
\begin{proof}
By Theorem 7 (1) in \cite{CoxPerezSendra2020}, $\deg(\cS)^2=\mapdeg(\cS)+\mult(\myB(\cS))$. Since $\cS\in \sol(\cP)$, $\mapdeg(\cS)=\mapdeg(\cP)$ and hence the statement follows.
\end{proof}

Observe that the minimal expected degree of an element in $\sol(\cP)$ is $\lceil\sqrt{\mapdeg(\cP)}\rceil$ which holds when $\myB(\cS)=\emptyset$.
In the next lemma, we analyze the case where the birational parametrization $\cQ$ has empty base locus.

\para

\begin{lemma}\label{lem:grado-sin-puntos-base}
Let us assume that $\myS$ admits a birational parametrization without base points.  Then, there exists
 $\cS\in \sol(\cP)$ such that $$\deg(\cS)=\dfrac{\deg(\cP)}{\sqrt{\deg(\myS)}}.$$
\end{lemma}
\begin{proof}
Let $\cQ$ be a birational parametrization of $\myS$ such that $\myB(\cQ)=\emptyset$. Then, $(\cQ,\cQ^{-1}\circ \cP)$ solves the
birational reparametrization problem. So, by Theorem \ref{thrm:caracterizacion}, $\cS:=\cQ^{-1}\circ \cP\in \sol(\cP)$. Moreover,
by Corollary 10 in \cite{CoxPerezSendra2020}, $\deg(\cP)=\deg(\cQ)\deg(\cS)$. Now, by Theorem 3 in \cite{CoxPerezSendra2020} applied to $\cQ$, $\deg(\cQ)=\sqrt{\deg(\myS)}$. This concludes the proof.
\end{proof}

\subsection{General algorithm}\label{subsec:generalAlgoritm}
In this subsection, combining all previous results and ideas,  the  general algorithm is derived.   For this purpose, we  identify the set of all projective curves, including multiple component curves, of a fixed degree $d$, with the projective space (see \cite{Miranda},  \cite{SWP} or  \cite{Walker} for further details)
\begin{equation}\label{eq-Vd}
\myV_d:=\mathbb{P}^{\frac{d (d+3)}{2}}(\K).
\end{equation}
More precisely, we identify the projective curves of degree $d$ with the forms in $\K[\ot]$ of degree $d$, up to multiplication by non-zero $\K$-elements. Now, these forms are identified with the elements in $\myV_d$ corresponding to their coefficients, after fixing an order of the monomials. By abuse of notation, we will refer to the elements in $\myV_d$ by either their tuple of coefficients, or the associated form, or the corresponding curve.   In the sequel, let us denote by
\begin{equation}\label{eq-Egeneral}
E_d(\Lambda,\ot)\in \K[\lambda,\ot]
\end{equation}
the generic $\ot$--homogeneous polynomial in $\myV_d$, where $\Lambda$ is the tuple of undetermined coefficients.

\para

First, we outline the auxiliary   Algorithm \ref{alg:subSol1}  that computes the set $\cU$ of curves, of a fixed given degree $d$, satisfying Theorem \ref{thrm:caracterizacion} (2a). In the sequel, we describe the basic ideas of our approach.
We will see that $\cU$ is an open subset, maybe empty, of a closed subset $\myW$  of $\myV_d\times \myV_d \times \myV_d$.  For this purpose, let
\begin{equation}\label{eq-Edi}
E_{d}^{i}:=E(\Lambda_{i},\ot), \,\, i\in \{1,2,3\}
\end{equation}
where $\Lambda_{1},\Lambda_2,\Lambda_{3}$ are different tuples of undetermined coefficients; we will use the notation $\overline{\Lambda}:=
(\overline{\Lambda}_1,\overline{\Lambda}_2,\overline{\Lambda}_3)$. We will find  $U_1,\ldots,U_r\in \K[\Lambda_1,\Lambda_2,\Lambda_3]$ and a variety $\myW\subset \myV_{d}^{3}$ such that   $\cU=\myW\setminus \cup_{i=1}^{r}\mathbb{V}(U_i)$.
The reasoning is as follows.     $\cE:=(E_{d}^{1}:E_{d}^{2}:E_{d}^{3})$ defines a generic rational map from $\proj2$ into $\proj2$. We look for those $\cE$ such that $\myF_g(\cE)=\myF_g(\cP)$.  To achieve this, we will ensure that $\myF_g(\cP)\subset \myF_g(\cE)$ and that $\mapdeg(\cP)=\mapdeg(\cE)$.

\para

First let $U_1$ contain the parameters in $\Lambda_3$; this   guarantees that for  $\cE\in \myV_{d}^{3}\setminus \mathbb{V}(U_1)$, $E_{d}^{3}$ is not zero.  We need to impose first that $\myF_g(\cP)\subset \myF_g(\cE)$. But, since $E_{d}^{3}\neq 0$, according to Subsection \ref{subsec:generic-fibra}, it is enough to ask that  $\myF_g(\cP_a)\subset \myF_g(\cE_a)$, where $e_{d}^{i}:=E_{d}^{i}(\Lambda_{i},t_1,t_2,1)$ and
\begin{equation}\label{eq:Ea}
\cE_a:=\left(\dfrac{e_{d}^{1}}{e_{d}^{3}},
\dfrac{e_{d}^{2}}{e_{d}^{3}}\right):\mathbb{A}^{2}(\K)\dashrightarrow \mathbb{A}^{2}(\K).
\end{equation}
Let $U_2$ be the set containing  the conditions on the parameters $\overline{\Lambda}$ that ensure  that the rank of the Jacobian matrix of $\cE_a$ is smaller than 2; this guarantees that for  $\cE\in \myV_{d}^{3}\setminus \cup_{i=1}^{2} \mathbb{V}(U_i)$ the hypothesis in Subsection 1.2. in \cite{PS-grado}  on the linear independency is satisfied (see Remark \ref{rem:hypothses}) and that $\cE_a$ is dominant. We observe that the denominator of the determinant Jacobian is a power of $e_{d}^{3}$ and, because of $U_1$, the determinant of the Jacobian is always well-defined. Thus, in $U_2$, we only need to collect the coefficient w.r.t. $\{t_1,t_2\}$ of the  {numerator of the} determinant of the Jacobian.
Now, we introduce the polynomials  (compare to \eqref{eq-GS})
\begin{equation}\label{eq:GE}
\tilde{G}_{i}^{\cE}:=e_{d}^{i}(\Lambda_{i},h_1,h_2)e_{d}^{3}(\Lambda_{3},t_1,t_2)- e_{d}^{i}(\Lambda_{i},t_1,t_2)e_{d}^{3}(\Lambda_{3},h_1,h_2),\,\,i\in \{1,2\},
\end{equation}
 and we require that $\tilde{G}_{1}^{\cE}(\Lambda_{1},\Lambda_{3},h_1,h_2,\aa)=0$ for every $\aa\in \myF_g(\cP_a)$. We observe that
 $(h_1,h_2)\in \myF_g(\cP_a)$ and that $\tilde{G}_{1}^{\cE}(\Lambda_{1},\Lambda_{3},h_1,h_2,h_1,h_2)=0$.
 So, we may work with $\myF_g(\cP_a)\setminus \{(h_1,h_2)\}$. Furthermore,  $\myF_g(\cP_a)=\mathbb{V}(\myG)$, see \eqref{eq-F-BGBuena}, where $\myG$ is as in \eqref{eq-GBbuena}. Therefore, let us replace  the univariate polynomial $u(t_1)$ of  $\myG$, that is the square-free part of $R_1$ (see \eqref{eq-R1R2}),  by $$u^{*}(t_1):=u(t_1)/(t_1-h_1)$$ in $\myG$; let
 \begin{equation}\label{eq-GBstar}
 \myGG:=\{u^*(t_1),t_2-v(t_1)\}
 \end{equation} be the resulting set. Then,  $\mathscr{G}^*$ is still a Gr\"obner basis and
 \begin{equation}\label{eq-myFstar}
 \myF_g(\cP_a)\setminus \{(h_1,h_2)\}=\mathbb{V}(\mathscr{G}^*).
 \end{equation}
Thus, to achieve the condition above, one has to require
 the normal form $N_i$ of $\tilde{G}_{i}^{\cE}, i\in \{1,2\}$, w.r.t. $\mathscr{G}^*$ to be zero. We have the following result.

\para

\begin{lemma}\label{lemma-NF}
 Let $N_i$ be  the normal form   of $\tilde{G}_{i}^{\cE}, i\in \{1,2\}$, w.r.t. $\mathscr{G}^*$. It holds that
 \begin{enumerate}
 \item $N_i, i\in \{1,2\}$, is the remainder of the division of $\tilde{G}_{i}^{\cE}(t_1,v(t_1))$ by $u^*(t_1)$ w.r.t. $t_1$.
 \item For $\cE\in \myV_{d}^{3}\setminus \cup_{i=1}^{2} \mathbb{V}(U_i)$, $N_1,N_2$ specializes properly.
 \end{enumerate}
 \end{lemma}
 \begin{proof}
 It follows from the special form of the Gr\"obner basis  $\myGG$, and from the fact that the polynomials in $\myGG$ does not depend on $\Lambda_i$
 \end{proof}

 \para

 If $ {N_1}$ is not zero and does not depend on $\Lambda_{1},\Lambda_{3}$, then there exists no element in $\cE\in \myV_{d}^{3}\setminus \cup_{i=1}^{2} \mathbb{V}(U_i)$ satisfying Theorem \ref{thrm:caracterizacion} (2) . So, let us assume that $ {N_1}$  does depend on
$\Lambda_{1},\Lambda_3$.  The condition $N_1=0$ provides a finite set  $C_{1,3}(\Lambda_1,\Lambda_3)$  of polynomials in the parameters $\Lambda_{1},\Lambda_{3}$, namely the non-zero coefficients w.r.t. $\{h_1,h_2\}$ of the coefficients w.r.t. $\{t_1,t_2\}$ of $N_1$. Let  $C_{2,3}(\Lambda_2,\Lambda_3)$ the corresponding set to $N_2$. By symmetry $C_{2,3}(\Lambda_2,\Lambda_3)= C_{1,3}(\Lambda_2,\Lambda_3)$.

As a consequence, we have the following result.

\para

\begin{lemma}
Let $C_{1,3}$ and $U_i$ be as above. For every $(E_1:E_3),(E_2:E_3)\in  \mathbb{V}(C_{1,3})$
such that $\cE:=(E_1:E_2:E_3)\not\in  \cup_{i=1}^{2} \mathbb{V}(U_i)$ it holds that
$\myF_g(\cP)\subset \myF_g(\cE)$.
\end{lemma}

 \para

We observe that the polynomials in $C_{1,2}$ are easy to handle.

 \para

 \begin{lemma}
The polynomials in $C_{1,3}$ are bilinear  forms in the undetermined parameters $\Lambda_{1}$ and $\Lambda_{3}$.
 \end{lemma}
 \begin{proof}
 The result follows taking into account the reduction process for computing the normal form, that  $G_{1}^{\cE}$ is  bilinear in the parameters $\Lambda_{1},\Lambda_{3}$, and that no polynomial in $\myGG$ depends on $\Lambda_{1},\Lambda_{3}$.
 \end{proof}

 \para

 In this situation, solving the algebraic system of bilinear forms $\{h(\Lambda_{1},\Lambda_{3})=0\}_{h\in C_{1,3}}$ yields in general to a finite set of  (parametric) solutions which provides a description, by means of rational generic elements, on the irreducible components of
 $\mathbb{V}(C_{1,3})$. In this case, the denominators of each generic point  will be included in the set $U_2$; or equivalently the condition requiring that the projective point does not collapse to the zero tuple.  Now, for every generic solution $(E_1,E_3)$  of $\mathbb{V}(C_{1,3})$, let $\ov=(E_1:E_2:E_3)$ where $E_2$ is the result of substituting in $E_1$ the parameters $\Lambda_1$ by $\Lambda_2$.
 We denote by $\cV_{1,2,3}$ be the set of all points constructed in this way.

 \para

 Let $\ov\in \cV_{1,2,3}$.
 We look for the conditions such that $\mapdeg(\cP)=\mapdeg(\cE)$. For this we use Subsection \ref{subsec:fiber-S}. More precisely,
 let (see \eqref{eq:GE})
 \begin{equation}\label{eq:GEbreve}
\breve{G}_{i}^{\cE}:=\tilde{G}_{i}^{\cE}(\ov,h_1,h_2,t_1,t_2),\,\,i\in \{1,2\},
\end{equation}
and
 \begin{equation}\label{eq-R1Ebreve}
 \breve{R}_{1}^{\cE}=\mathrm{PrimPart}_{\{t_1,t_2\}}(\mathrm{PrimPart}_{\{h_1,h_2\}}(\Res_{t_2}(\breve{G}_{1}^{\cE},\breve{G}_{2}^{\cE}))).
 \end{equation}
 Let $\breve{R}_{1}^{\cE}$ be expressed as
 \begin{equation}\label{eq-R1Eb}
 \breve{R}_{1}^{\cE}=\sum c_{i}(\Lambda_1,\Lambda_{2},\Lambda_{3},h_1,h_2) t^i.
 \end{equation}
 Since $\mapdeg(\cE)=\deg_{t_1}(\breve{R}_{1}^{\cE})$, let $C^*$ be the set of all non-zero coefficients of $c_{i}$ w.r.t. $\{h_1,h_2\}$ and $i>\mapdeg(\cP)$. In addition, let $U_3$ be the set of all non-zero coefficients of $c_{\mapdeg(\cP)}$ w.r.t. $\{h_1,h_2\}$.  Then, we have the following result.

 \para

\begin{theorem}\label{theorem:CU}
Let $C^*$ and $U_i$ be as above. For every $\cE\in \mathbb{V}(C^*)\setminus  \cup_{i=1}^{3} \mathbb{V}(U_i)$ it holds that
$\myF_g(\cP)= \myF_g(\cE)$.
\end{theorem}

 \para

   In the following algorithm, we summarize all these ideas.

 \vspace*{3mm}

\begin{breakablealgorithm}
\caption{\textsf{SolSpace}}
\label{alg:subSol1}
\begin{algorithmic}[1]
	\REQUIRE The set of polynomials $\myGG:=\{u^*(t_1),t_2-v(t_1)\}$ as in \eqref{eq-GBstar}, and a positive integer $d\geq \lceil \sqrt{\deg_{t_1}(R_1) }\rceil$.
	\ENSURE  \textbf{Either} A pair of sets  $[C,[U_1,U_2,U_3]]$   as in Theorem \ref{theorem:CU} \textbf{or} the empty set.
	\STATE 
	Set $C:=\emptyset$, $U_i=\emptyset$ for $i=1,2$.
	\STATE Let  $E_d(\Lambda,\ot)$ be as in \eqref{eq-Egeneral}. Set $E_{d}^{i}:=E_d(\Lambda_i,\ot), i\in \{1,2,3\}$, where $\Lambda_i$ is a tuple of undetermined coefficients.
	\STATE Include in $U_1$ the entries of $\Lambda_{3}$,   and in $U_2$
	the non-zero coefficients w.r.t. $\{t_1,t_2\}$ of the determinant of   the Jacobian of $\cE_a$ (see \eqref{eq:Ea}).
	\IF{$U_2$ contains a non-zero constant} \STATE Set $U_2=\{1\}$.
	\ENDIF
	\STATE Let $\tilde{G}_{i}^{\cE}$ be as in \eqref{eq:GE} for $i\in \{1,2\}$.
	\STATE Compute the normal form $N_1$ of $\tilde{G}_{1}^{\cE}$ w.r.t. $\myGG$ as in Lemma \ref{lemma-NF}.
	\IF{$N_1\in \K(h_1,h_2)[t_1,t_2]\setminus \{0\}$} \STATE \textbf{return} $\emptyset$
	\ENDIF
	\STATE\label{hola2} Let   $C_{1,3}$ be the set  non-zero coefficients w.r.t. $\{h_1,h_2\}$ of the coefficients w.r.t. $\{t_1,t_2\}$ of $N_1$.
	\STATE Solve the system of bilinear forms $C_{1,3}$. Let $\cV\subset \myV_{d}^{2}$ be the set of solutions.
	\STATE For each $(E_1,E_3)\in \cV$, let $\ov=(E_1:E_2:E_3)$ where $E_2$ is the result of substituting in $E_1$ the parameters $\Lambda_1$ by $\Lambda_2$.  Let $\cV_{1,2,3}$ be the set of all points constructed in this way.
	\STATE Delete from $\cV_{1,2,3}$ those solutions  vanishing all polynomials in, at least, one of the sets $U_i, i=1,2$. Let $\cV_{1,2,3}$ be the resulting set.
	\IF{$\cV_{1,2,3}=\emptyset$} \STATE \textbf{return} $\emptyset$ \ENDIF
	\STATE Set $C=\emptyset$ and $U_3=\emptyset$.
	\FOR{each irreducible component of $\cV_{1,2,3}$, choose $\ov$ a generic element and} \label{step-for-v}
	\STATE Compute $\breve{G}_{i}^{\cE}:=\tilde{G}_{i}^{\cE}(\ov,h_1,h_2,t_1,t_2)$ for $i=1,2$ and $\breve{R}_{1}^{\cE}$ as in \eqref{eq-R1Ebreve}.
	\IF{$\deg_{t_1}(\breve{R}_{1}^{\cE})<\deg_{t_1}(R_1)$} \STATE go to Step \ref{step-for-v}
	\ELSE
	\STATE Let $C^*$ be the set of all non-zero coefficients of $c_{i}$ in \eqref{eq-R1Eb}, w.r.t. $\{h_1,h_2\}$, and $i>\deg_{t_1}(R_1)$.
	\IF{an element in $C^*$ is a non-zero constant} \STATE go to Step \ref{step-for-v}
	\ELSE
	\STATE Replace $C$ by $C\cup C^*$.
	\STATE Let $U^*$ be  the set of all non-zero coefficients of $c_{\deg_{t_1}(R_1)}$ w.r.t. $\{h_1,h_2\}$.
	 \IF{$U^*$ contains a non-zero constant} \STATE Set $U_3=\{1\}$.  \ELSE
	\STATE Replace $U_3$ by $U_3\cup U^*$.
		\ENDIF
        \ENDIF
	\ENDIF
	\ENDFOR
	\STATE Solve   $C$. Let $\cV\subset \myV_{d}^{3}$ be the set of solutions.
	\STATE Delete from $\cV$ those solutions    vanishing all polynomials in, at least, one of the sets $U_i$. Let $\cV$ be the resulting set.
	\STATE \textbf{return} $[C,U:=[U_1,U_2,U_3]]$.		
\end{algorithmic}
\end{breakablealgorithm}

\para

Finally,   the general algorithm is derived. We assume that the algorithm in \cite{sonia} has been applied and has not provide an answer to the problem.

 \para

\begin{breakablealgorithm}
\caption{\textsf{General Algorithm}}
\label{alg:general}
\begin{algorithmic}[1]
	\REQUIRE $\cP$ a rational surface parametrization as in \eqref{eq-P}.
	\ENSURE   $\cS$ and $\cQ$ a solution to the birational reparametrization problem for $\cP$.
	\STATE Compute $\ell=\lceil\sqrt{\mapdeg(\cP)}\rceil$ (see \eqref{eq-degree} and Lemma \ref{lem:grado-general}).
	\STATE Compute $\myG^*$ as in \eqref{eq-GBstar}.
	\FOR{$d\geq \ell$}\label{pasoFor} \STATE\label{paso:subalg1}  Apply Algorithm \ref{alg:subSol1} with input $\myGG$ and $d$.
	\IF{the output of Algorithm \ref{alg:subSol1}  is the empty set}
	\STATE set $d:=d+1$ and go to Step \ref{pasoFor}
	\ELSE
         \STATE Let $[C,[U_1,U_2,U_3]]$ be the output of Algorithm \ref{alg:subSol1}.
         \STATE Find $\ov\in \mathbb{V}(C) \setminus \cup_{i=1}^{3} \mathbb{V}(U_i)$.
          \STATE  $\cS:=(E_{d}^{1}(\ov):E_{d}^{2}(\ov):E_d^{3}(\ov))$ see   \eqref{eq-Edi}.
           \FOR{$j\in \{1,2,3\}$}
         \STATE Determine $\Phi_j $ see \eqref{eq-Phi}.
         \STATE Compute the implicit equation $H_j:=A_{j,1}(x,y)-A_{j,0}(x,y)z$ of the surface parametrization $\Phi_j$.
         \ENDFOR
         \STATE Set $\cQ:=(A_{1,1}/A_{1,0}, A_{2,1}/A_{2,0},A_{3,1}/A_{3,0})$.
	 \ENDIF
	 \STATE \textbf{return} $[\cS,\cQ]$.
	\ENDFOR
	\end{algorithmic}
\end{breakablealgorithm}

 \para

We finish this section with an example where we illustrate how the general algorithm works.

 \para

\begin{example}\label{example:general}
Let $\cP(\ot)= (t_1^3+t_2t_{3}^{2}: t_1^3: t_2t_{3}^{2}:t_3^3)$ be the input of Algorithm \ref{alg:general}; that is, it is a  rational projective parametrization of an algebraic surface $\myS$ as in \eqref{eq-P}. The parametrization is quite simple and   one could implicitize  and  compute a proper rational parametrization. Nevertheless, let us use this parametrization to illustrate carefully Algorithm \ref{alg:general} and Algorithm \ref{alg:subSol1}. For this purpose, in Steps 1 and 2 of Algorithm \ref{alg:general},  we first compute $\ell=\sqrt{\mapdeg(\cP)}=\sqrt{3}$ and the set of polynomials $$\myGG:=\{ h_{1}^{2}+h_{1} t_{1}+t_{1}^{2}, t_2-h_2\}$$ as in \eqref{eq-GBstar}. Thus, $u^*(t_1)= h_{1}^{2}+h_{1} t_{1}+t_{1}^{2}$ and $v(t_1)=h_2$.

The for-loop in Step 3  starts with $d=2$. We apply   Algorithm \ref{alg:subSol1} to $\myGG$ and $d=2$:

\para

\noindent \textsf{Algorithm  \ref{alg:subSol1}} starts.

\para

\noindent \textsf{Step 1.}  $C=U_1=U_2=\emptyset$.
\vspace*{1mm}

\noindent  \textsf{Step 2.}  $E_2^i(\Lambda_i, t_1, t_2):= \lambda_{i,1}t_1^2+\lambda_{i,2}t_2^2+\lambda_{i,3}t_1t_2+\lambda_{i,4}t_1t_3+\lambda_{i,5}t_2t_3+\lambda_{i,6}t_3^2,\,i\in \{1,2,3\}$.

 \vspace*{1mm}

\noindent \textsf{Steps 3-6.}  $U_1=\{\lambda_{3,1},\ldots,\lambda_{3,6}\}$. For determining $U_2$, let
    \[\cE_a=\left(\frac{e_2^1}{e_2^3},\frac{e_2^2}{e_2^3}\right),\]
 where $e_2^i=E_2^i(t_1, t_2, 1)$ for $i=1,2,3$. Let   $T$ be the determinant of the Jacobian of $\cE_a$ w.r.t. $\{t_1,t_2\}$. Then, $U_2$ collects the coefficients of the numerator of $T$ w.r.t. $\{t_1,t_2\}$. In this case, $U_2$ consists in ten 3-linear polynomials in
 $\lambda_{1,i}, \lambda_{2,j}, \lambda_{3,k}$.
   \vspace*{1mm}

\noindent \textsf{Step
  7.} Let  $\tilde{G}_i^{\cE}=e_2^i(\Lambda_i, h_1, h_2)e_2^3(\Lambda_3, t_1, t_2)-e_2^3(\Lambda_3, h_1, h_2)e_2^i(\Lambda_i, t_1, t_2),\,i=1,2,3.$

   \vspace*{1mm}

\noindent \textsf{Step 8.} We get that the normal form $N_1$of $\tilde{G}_1^{\cal E}$ w.r.t ${\cal G}^{\star}$ is
   \vspace*{1.5mm}

   \noindent
   $N_1=-h_{1}^{3} h_{2} \lambda_{1,1} \lambda_{3,3}+h_{1}^{3} h_{2} \lambda_{1,3} \lambda_{3,1}-2 h_{1}^{2} h_{2}^{2} \lambda_{1,1} \lambda_{3,2}+2 h_{1}^{2} h_{2}^{2} \lambda_{1,2} \lambda_{3,1}-2 h_{1}^{2} h_{2} \lambda_{1,1} \lambda_{3,3} t_{1}+2 h_{1}^{2} h_{2} \lambda_{1,3} \lambda_{3,1} t_{1}+h_{1} h_{2}^{3} \lambda_{1,2} \lambda_{3,3}-h_{1} h_{2}^{3} \lambda_{1,3} \lambda_{3,2}-h_{1} h_{2}^{2} \lambda_{1,1} \lambda_{3,2} t_{1}+h_{1} h_{2}^{2} \lambda_{1,2} \lambda_{3,1} t_{1}-h_{2}^{3} \lambda_{1,2} \lambda_{3,3} t_{1}+h_{2}^{3} \lambda_{1,3} \lambda_{3,2} t_{1}-h_{1}^{3} \lambda_{1,1} \lambda_{3,4}+h_{1}^{3} \lambda_{1,4} \lambda_{3,1}-2 h_{1}^{2} h_{2} \lambda_{1,1} \lambda_{3,5}+2 h_{1}^{2} h_{2} \lambda_{1,5} \lambda_{3,1}-2 h_{1}^{2} \lambda_{1,1} \lambda_{3,4} t_{1}+2 h_{1}^{2} \lambda_{1,4} \lambda_{3,1} t_{1}+h_{1} h_{2}^{2} \lambda_{1,2} \lambda_{3,4}-h_{1} h_{2}^{2} \lambda_{1,3} \lambda_{3,5}-h_{1} h_{2}^{2} \lambda_{1,4} \lambda_{3,2}+h_{1} h_{2}^{2} \lambda_{1,5} \lambda_{3,3}-h_{1} h_{2} \lambda_{1,1} \lambda_{3,5} t_{1}+h_{1} h_{2} \lambda_{1,5} \lambda_{3,1} t_{1}-h_{2}^{2} \lambda_{1,2} \lambda_{3,4} t_{1}+h_{2}^{2} \lambda_{1,3} \lambda_{3,5} t_{1}+h_{2}^{2} \lambda_{1,4} \lambda_{3,2} t_{1}-h_{2}^{2} \lambda_{1,5} \lambda_{3,3} t_{1}-2 h_{1}^{2} \lambda_{1,1} \lambda_{3,6}+2 h_{1}^{2} \lambda_{1,6} \lambda_{3,1}-h_{1} h_{2} \lambda_{1,3} \lambda_{3,6}-h_{1} h_{2} \lambda_{1,4} \lambda_{3,5}+h_{1} h_{2} \lambda_{1,5} \lambda_{3,4}+h_{1} h_{2} \lambda_{1,6} \lambda_{3,3}-h_{1} \lambda_{1,1} \lambda_{3,6} t_{1}+h_{1} \lambda_{1,6} \lambda_{3,1} t_{1}+h_{2} \lambda_{1,3} \lambda_{3,6} t_{1}+h_{2} \lambda_{1,4} \lambda_{3,5} t_{1}-h_{2} \lambda_{1,5} \lambda_{3,4} t_{1}-h_{2} \lambda_{1,6} \lambda_{3,3} t_{1}-h_{1} \lambda_{1,4} \lambda_{3,6}+h_{1} \lambda_{1,6} \lambda_{3,4}+\lambda_{1,4} \lambda_{3,6} t_{1}-\lambda_{1,6} \lambda_{3,4} t_{1}
$

  \vspace*{1.5mm}

\noindent \textsf{Step 9-11.}   Since $N_1 \not\in {\Bbb K}(h_1, h_2)[t_1, t_2]\setminus\{0\}$, we continue in Step 12.   \vspace*{1mm}

\noindent \textsf{Step  12.} We have that
  \vspace*{1.5 mm}

\noindent
$C_{1,3}=\{-2 \lambda_{1,1} \lambda_{3,2}+2 \lambda_{1,2} \lambda_{3,1},
-\lambda_{1,1} \lambda_{3,2}+\lambda_{1,2} \lambda_{3,1},
-2 \lambda_{1,1} \lambda_{3,3}+2 \lambda_{1,3} \lambda_{3,1},
-\lambda_{1,1} \lambda_{3,3}+\lambda_{1,3} \lambda_{3,1},
-2 \lambda_{1,1} \lambda_{3,4}+2 \lambda_{1,4} \lambda_{3,1},
-\lambda_{1,1} \lambda_{3,4}+\lambda_{1,4} \lambda_{3,1},
-2 \lambda_{1,1} \lambda_{3,5}+2 \lambda_{1,5} \lambda_{3,1},
-\lambda_{1,1} \lambda_{3,5}+\lambda_{1,5} \lambda_{3,1},
-2 \lambda_{1,1} \lambda_{3,6}+2 \lambda_{1,6} \lambda_{3,1},
-\lambda_{1,1} \lambda_{3,6}+\lambda_{1,6} \lambda_{3,1},
-\lambda_{1,2} \lambda_{3,3}+\lambda_{1,3} \lambda_{3,2},
\lambda_{1,2} \lambda_{3,3}-\lambda_{1,3} \lambda_{3,2},
-\lambda_{1,4} \lambda_{3,6}+\lambda_{1,6} \lambda_{3,4},
\lambda_{1,4} \lambda_{3,6}-\lambda_{1,6} \lambda_{3,4},
-\lambda_{1,2} \lambda_{3,4}+\lambda_{1,3} \lambda_{3,5}+\lambda_{1,4} \lambda_{3,2}-\lambda_{1,5} \lambda_{3,3}
,
\lambda_{1,2} \lambda_{3,4}-\lambda_{1,3} \lambda_{3,5}-\lambda_{1,4} \lambda_{3,2}+\lambda_{1,5} \lambda_{3,3}
,
-\lambda_{1,3} \lambda_{3,6}-\lambda_{1,4} \lambda_{3,5}+\lambda_{1,5} \lambda_{3,4}+\lambda_{1,6} \lambda_{3,3}
,
\lambda_{1,3} \lambda_{3,6}+\lambda_{1,4} \lambda_{3,5}-\lambda_{1,5} \lambda_{3,4}-\lambda_{1,6} \lambda_{3,3}
\}$

\vspace*{1.5mm}

  \noindent
  is the set non-zero coefficients w.r.t. $\{h_1, h_2\}$ of the coefficients w.r.t.
$\{t_1, t_2\}$  of $N_1$.

 \vspace*{1mm}

\noindent \textsf{Step  13.} We solve the system of bilinear forms $C_{1,3}$. We get $\cV$ as the set solutions
     \vspace*{1mm}

  \noindent $  \{\{\lambda_{1,1} = 0, \lambda_{1,2} = \lambda_{1,2},
\lambda_{1,3} = 0, \lambda_{1,4} = 0, \lambda_{1,5} = \lambda_{1,5},
\lambda_{1,6} = \lambda_{1,6}, \lambda_{3,1} = 0, \lambda_{3,2} =
\lambda_{3,2}, \lambda_{3,3} = 0, \lambda_{3,4} = 0, \lambda_{3,5} =
\lambda_{3,5}, \lambda_{3,6} = \lambda_{3,6}\}, \{\lambda_{1,1}
 = 0, \lambda_{1,2} = \lambda_{1,2}, \lambda_{1,3} = \lambda_{1,3},
\lambda_{1,4} = \lambda_{1,4}, \lambda_{1,5} = \lambda_{1,5},
\lambda_{1,6} =
- {\lambda_{1,4} \left(\lambda_{1,2} \lambda_{1,4}-\lambda_{1,3} \lambda_{1,5}\right)}/{\lambda_{1,3}^{2}}
, \lambda_{3,1} = 0, \lambda_{3,2} = 0, \lambda_{3,3} = 0,
\lambda_{3,4} = \lambda_{3,4}, \lambda_{3,5} =
{\lambda_{1,2} \lambda_{3,4}}/{\lambda_{1,3}}, \lambda_{3,6} =
- {\lambda_{3,4} \left(\lambda_{1,2} \lambda_{1,4}-\lambda_{1,3} \lambda_{1,5}\right)}/{\lambda_{1,3}^{2}}
\},  \{\lambda_{1,1} = 0, \lambda_{1,2} =
 {\lambda_{1,3} \lambda_{3,2}}/{\lambda_{3,3}}, \lambda_{1,3} =
\lambda_{1,3}, \lambda_{1,4} = \lambda_{1,4}, \lambda_{1,5} =
- {\lambda_{1,3} \lambda_{3,2} \lambda_{3,4}-\lambda_{1,3} \lambda_{3,3} \lambda_{3,5}-\lambda_{1,4} \lambda_{3,2} \lambda_{3,3}}/{\lambda_{3,3}^{2}}
, \lambda_{1,6} =
- {\lambda_{1,4} \left(\lambda_{3,2} \lambda_{3,4}-\lambda_{3,3} \lambda_{3,5}\right)}/{\lambda_{3,3}^{2}}
, \lambda_{3,1} = 0, \lambda_{3,2} = \lambda_{3,2}, \lambda_{3,3} =
\lambda_{3,3}, \lambda_{3,4} = \lambda_{3,4}, \lambda_{3,5} =
\lambda_{3,5}, \lambda_{3,6} =
- {\lambda_{3,4} \left(\lambda_{3,2} \lambda_{3,4}-\lambda_{3,3} \lambda_{3,5}\right)}/{\lambda_{3,3}^{2}}
 \}, \{\lambda_{1,1} = \lambda_{1,1}, \lambda_{1,2} =
\lambda_{1,2}, \lambda_{1,3} = 0, \lambda_{1,4} = 0, \lambda_{1,5} =
\lambda_{1,5}, \lambda_{1,6} = \lambda_{1,6}, \lambda_{3,1} =
\lambda_{3,1}, \lambda_{3,2} =
 {\lambda_{1,2} \lambda_{3,1}}/{\lambda_{1,1}}, \lambda_{3,3} = 0,
\lambda_{3,4} = 0, \lambda_{3,5} =
 {\lambda_{1,5} \lambda_{3,1}}/{\lambda_{1,1}}, \lambda_{3,6} =
 {\lambda_{1,6} \lambda_{3,1}}/{\lambda_{1,1}}\}, \{
\lambda_{1,1} = \lambda_{1,1}, \lambda_{1,2} = \lambda_{1,2},
\lambda_{1,3} = \lambda_{1,3}, \lambda_{1,4} = \lambda_{1,4},
\lambda_{1,5} = \lambda_{1,5}, \lambda_{1,6} = \lambda_{1,6},
\lambda_{3,1} = 0, \lambda_{3,2} = 0, \lambda_{3,3} = 0, \lambda_{3,4}
 = 0, \lambda_{3,5} = 0, \lambda_{3,6} = 0\},  \{\lambda_{1,1} =
 {\lambda_{1,3} \lambda_{3,1}}/{\lambda_{3,3}}, \lambda_{1,2} =
 {\lambda_{1,3} \lambda_{3,2}}/{\lambda_{3,3}}, \lambda_{1,3} =
\lambda_{1,3}, \lambda_{1,4} =
 {\lambda_{1,3} \lambda_{3,4}}/{\lambda_{3,3}}, \lambda_{1,5} =
 {\lambda_{1,3} \lambda_{3,5}}/{\lambda_{3,3}}, \lambda_{1,6} =
 {\lambda_{1,3} \lambda_{3,6}}/{\lambda_{3,3}}, \lambda_{3,1} =
\lambda_{3,1}, \lambda_{3,2} = \lambda_{3,2}, \lambda_{3,3} =
\lambda_{3,3}, \lambda_{3,4} = \lambda_{3,4}, \lambda_{3,5} =
\lambda_{3,5}, \lambda_{3,6} = \lambda_{3,6} \}, \{
\lambda_{1,1} =  {\lambda_{1,4} \lambda_{3,1}}/{\lambda_{3,4}},
\lambda_{1,2} = {\lambda_{1,4} \lambda_{3,2}}/{\lambda_{3,4}},
\lambda_{1,3} = 0, \lambda_{1,4} = \lambda_{1,4}, \lambda_{1,5} =
 {\lambda_{1,4} \lambda_{3,5}}/{\lambda_{3,4}}, \lambda_{1,6} =
{\lambda_{1,4} \lambda_{3,6}}/{\lambda_{3,4}}, \lambda_{3,1} =
\lambda_{3,1}, \lambda_{3,2} = \lambda_{3,2}, \lambda_{3,3} = 0,
\lambda_{3,4} = \lambda_{3,4}, \lambda_{3,5} = \lambda_{3,5},
\lambda_{3,6} = \lambda_{3,6} \} \}.
$

    \vspace*{1.5mm}

\noindent $\cV\subset \myV_d\times \myV_d$ contains 7 generic solutions.

  \vspace*{1mm}

\noindent \textsf{Step 14}   Replicating the $\Lambda_i$-elements  in $\cV$ we construct $\cV_{1,2,3}$.

\vspace*{1mm}

\noindent \textsf{Step 15-18}.  We delete   from $\cV_{1,2,3}$  those solutions  vanishing all polynomials in, at least, one of the sets $U_i,\, i =
1,2$. We get that the resulting set is an empty set. So, the output of Algorithm \ref{alg:subSol1} is the empty set. Therefore,  we set $d=3$, and we go to Step 3 in Algorithm \ref{alg:subSol1} to repeat the above steps.

\para

\noindent \textsf{Algorithm  \ref{alg:subSol1}} starts.

\para

\noindent \textsf{Step 1.} This step was already performed in the first iteration of Algorithm  \ref{alg:subSol1}.

\vspace*{1mm}

\noindent \textsf{Step 2.} Let $E_3^i:=
 \lambda_{i,1} t_{3}  t_{1}^{2}+ \lambda_{i,2}t_{2}^{3}+ \lambda_{i,3}t_{2}^{2} t_{3}+\lambda_{i,4}t_{2} t_{3}  t_{1}+ \lambda_{i,5} t_{3}^{2}t_{1}+\lambda_{i,6}t_{2} t_{3}^{2} +\lambda_{i,7}t_{3}^{3} +\lambda_{i,8} t_{1}^{3}+t_{2} \lambda_{i,9} t_{1}^{2}+\lambda_{i,10}t_{2}^{2}  t_{1}
$, $i\in \{1,2,3\}$.
\vspace*{1mm}

\noindent \textsf{Step 3-6.} $U_1=\{\lambda_{3,1},\ldots,\lambda_{3, 10}\}$ and reasoning as above $U_2$ consists in 28 3-linear polynomials in $\lambda_{1,i},\lambda_{2,j},\lambda_{3,k}$.

\vspace*{1mm}

\noindent \textsf{Step 7.} Let
  $\tilde{G}_i^{\cal E}=e_3^i(\Lambda_i, h_1, h_2)e_3^3(\Lambda_3, t_1, t_2)-e_3^3(\Lambda_3, h_1, h_2)e_3^i(\Lambda_i, t_1, t_2),\,i=1,2,3.$

\vspace*{1mm}

\noindent \textsf{Step 8-11.} We compute the normal form $N_1$ of $\tilde{G}_1^{\cal E}$ w.r.t ${\cal G}^{\star}$ and  we check that
   $N_1 \not\in {\Bbb K}(h_1, h_2)[t_1, t_2]\setminus\{0\}$.
\vspace*{1mm}

\noindent \textsf{Step 12.}  We have that
\vspace*{1.5 mm}

\noindent $C_{1,3}=\{-2 \lambda_{1,1} \lambda_{3,5}+2 \lambda_{1,5} \lambda_{3,1},
-\lambda_{1,1} \lambda_{3,5}+\lambda_{1,5} \lambda_{3,1},
-2 \lambda_{1,1} \lambda_{3,7}+2 \lambda_{1,7} \lambda_{3,1},
-\lambda_{1,1} \lambda_{3,7}+\lambda_{1,7} \lambda_{3,1},
-2 \lambda_{1,1} \lambda_{3,8}+2 \lambda_{1,8} \lambda_{3,1},
-\lambda_{1,1} \lambda_{3,8}+\lambda_{1,8} \lambda_{3,1},
\lambda_{1,2} \lambda_{3,9}-\lambda_{1,9} \lambda_{3,2},
2 \lambda_{1,2} \lambda_{3,9}-2 \lambda_{1,9} \lambda_{3,2},
-\lambda_{1,2} \lambda_{3,10}+\lambda_{1,10} \lambda_{3,2},
\lambda_{1,2} \lambda_{3,10}-\lambda_{1,10} \lambda_{3,2},
-\lambda_{1,4} \lambda_{3,8}+\lambda_{1,8} \lambda_{3,4},
\lambda_{1,4} \lambda_{3,8}-\lambda_{1,8} \lambda_{3,4},
-\lambda_{1,5} \lambda_{3,7}+\lambda_{1,7} \lambda_{3,5},
\lambda_{1,5} \lambda_{3,7}-\lambda_{1,7} \lambda_{3,5},
-\lambda_{1,5} \lambda_{3,8}+\lambda_{1,8} \lambda_{3,5},
\lambda_{1,5} \lambda_{3,8}-\lambda_{1,8} \lambda_{3,5},
\lambda_{1,8} \lambda_{3,9}-\lambda_{1,9} \lambda_{3,8},
2 \lambda_{1,8} \lambda_{3,9}-2 \lambda_{1,9} \lambda_{3,8},
-\lambda_{1,8} \lambda_{3,10}+\lambda_{1,10} \lambda_{3,8},
\lambda_{1,8} \lambda_{3,10}-\lambda_{1,10} \lambda_{3,8},
-2 \lambda_{1,9} \lambda_{3,10}+2 \lambda_{1,10} \lambda_{3,9},
-\lambda_{1,9} \lambda_{3,10}+\lambda_{1,10} \lambda_{3,9},
-2 \lambda_{1,1} \lambda_{3,2}+2 \lambda_{1,2} \lambda_{3,1}+2 \lambda_{1,3} \lambda_{3,9}-2 \lambda_{1,9} \lambda_{3,3}
,
-\lambda_{1,1} \lambda_{3,2}+\lambda_{1,2} \lambda_{3,1}+\lambda_{1,3} \lambda_{3,9}-\lambda_{1,9} \lambda_{3,3}
,
-2 \lambda_{1,1} \lambda_{3,3}+2 \lambda_{1,3} \lambda_{3,1}+2 \lambda_{1,6} \lambda_{3,9}-2 \lambda_{1,9} \lambda_{3,6}
,
-\lambda_{1,1} \lambda_{3,3}+\lambda_{1,3} \lambda_{3,1}+\lambda_{1,6} \lambda_{3,9}-\lambda_{1,9} \lambda_{3,6}
,
-2 \lambda_{1,1} \lambda_{3,4}+2 \lambda_{1,4} \lambda_{3,1}+2 \lambda_{1,5} \lambda_{3,9}-2 \lambda_{1,9} \lambda_{3,5}
,
-\lambda_{1,1} \lambda_{3,4}+\lambda_{1,4} \lambda_{3,1}+\lambda_{1,5} \lambda_{3,9}-\lambda_{1,9} \lambda_{3,5}
,
-2 \lambda_{1,1} \lambda_{3,10}+2 \lambda_{1,4} \lambda_{3,9}-2 \lambda_{1,9} \lambda_{3,4}+2 \lambda_{1,10} \lambda_{3,1}
,
-\lambda_{1,1} \lambda_{3,10}+\lambda_{1,4} \lambda_{3,9}-\lambda_{1,9} \lambda_{3,4}+\lambda_{1,10} \lambda_{3,1}
,
-2 \lambda_{1,1} \lambda_{3,6}+2 \lambda_{1,6} \lambda_{3,1}+2 \lambda_{1,7} \lambda_{3,9}-2 \lambda_{1,9} \lambda_{3,7}
,
-\lambda_{1,1} \lambda_{3,6}+\lambda_{1,6} \lambda_{3,1}+\lambda_{1,7} \lambda_{3,9}-\lambda_{1,9} \lambda_{3,7}
,
-\lambda_{1,2} \lambda_{3,4}-\lambda_{1,3} \lambda_{3,10}+\lambda_{1,4} \lambda_{3,2}+\lambda_{1,10} \lambda_{3,3}
,
\lambda_{1,2} \lambda_{3,4}+\lambda_{1,3} \lambda_{3,10}-\lambda_{1,4} \lambda_{3,2}-\lambda_{1,10} \lambda_{3,3}
,
-\lambda_{1,4} \lambda_{3,7}-\lambda_{1,5} \lambda_{3,6}+\lambda_{1,6} \lambda_{3,5}+\lambda_{1,7} \lambda_{3,4}
,
\lambda_{1,4} \lambda_{3,7}+\lambda_{1,5} \lambda_{3,6}-\lambda_{1,6} \lambda_{3,5}-\lambda_{1,7} \lambda_{3,4}
,
-\lambda_{1,2} \lambda_{3,5}-\lambda_{1,3} \lambda_{3,4}+\lambda_{1,4} \lambda_{3,3}+\lambda_{1,5} \lambda_{3,2}-\lambda_{1,6} \lambda_{3,10}+\lambda_{1,10} \lambda_{3,6}
,
\lambda_{1,2} \lambda_{3,5}+\lambda_{1,3} \lambda_{3,4}-\lambda_{1,4} \lambda_{3,3}-\lambda_{1,5} \lambda_{3,2}+\lambda_{1,6} \lambda_{3,10}-\lambda_{1,10} \lambda_{3,6}
,
-\lambda_{1,3} \lambda_{3,5}+\lambda_{1,4} \lambda_{3,6}+\lambda_{1,5} \lambda_{3,3}-\lambda_{1,6} \lambda_{3,4}-\lambda_{1,7} \lambda_{3,10}+\lambda_{1,10} \lambda_{3,7}
,
\lambda_{1,3} \lambda_{3,5}-\lambda_{1,4} \lambda_{3,6}-\lambda_{1,5} \lambda_{3,3}+\lambda_{1,6} \lambda_{3,4}+\lambda_{1,7} \lambda_{3,10}-\lambda_{1,10} \lambda_{3,7}
\}.
$

\vspace*{1.5mm}

\noindent \textsf{Step  13.} We solve the system of bilinear forms $C_{1,3}$. We get $\cV$ as the set solutions. It contains 21 generic solutions.

   \vspace*{1mm}

\noindent \textsf{Step 14}   Replicating the $\Lambda_i$-elements  in $\cV$ we construct $\cV_{1,2,3}$.

\vspace*{1mm}

\noindent \textsf{Step 15-18.}
We delete from  $\cV_{1,2,3}$  those solutions  vanishing all polynomials in, at least, one of the sets $U_i,\, i =
1,2$. We get that $\cV_{1,2,3}$ only contains one solution, namely,

\vspace*{1mm}

\noindent $\cV_{1,2,3}=\{\lambda_{1,1} = 0, \lambda_{1,2} = \lambda_{1,2}, \lambda_{1,3} =
\lambda_{1,3}, \lambda_{1,4} = 0, \lambda_{1,5} = 0, \lambda_{1,6} =
\lambda_{1,6}, \lambda_{1,7} = \lambda_{1,7}, \lambda_{1,8} =
\lambda_{1,8}, \lambda_{1,9} = 0, \lambda_{1,10} = 0, \lambda_{3,1} = 0
, \lambda_{3,2} = \lambda_{3,2}, \lambda_{3,3} = \lambda_{3,3},
\lambda_{3,4} = 0, \lambda_{3,5} = 0, \lambda_{3,6} = \lambda_{3,6},
\lambda_{3,7} = \lambda_{3,7}, \lambda_{3,8} = \lambda_{3,8},
\lambda_{3,9} = 0, \lambda_{3,10} = 0\}
$

\vspace*{1.5mm}

\noindent \textsf{Step 19.} We initialize $C:=\emptyset$ and $U_3:=\emptyset$.

\vspace*{1mm}

\noindent \textsf{Step 20.}  Since $\#(\cV_{1,2,3})=1$ the for-loop consists only in one iteration.

\vspace*{1mm}

\noindent \textsf{Steps 21-24.} We compute the specializations $\breve{G}_{1}^{\cE}, \breve{G}_{2}^{\cE}$ as well as the resultant $\breve{R}_{1}^{\cE}$. $\deg_{t_1}(\breve{R}_{1}^{\cE})=9$. So we go to Step 26.

\vspace*{1mm}

\noindent \textsf{Steps 25-29.} We compute $C^*$ (which is $C$ in this case) that consists in 21 polynomials

\vspace*{1mm}

\noindent \textsf{Steps 30-38.} $U_3$ consists in 18 polynomials.

\vspace*{1mm}

\noindent \textsf{Step 39.} We solve $C^*$ and we get $\cV$ with 12 generic solutions.

\vspace*{1mm}

\noindent \textsf{Steps 40.}  Filtering the solutions with $U_1,U_2,U_3$ we get  that  the new $\cV$ contains 3 solutions, namely

\vspace*{1.5 mm}

\noindent $\cV=\{\{\lambda_{1,1} = \lambda_{1,1}, \lambda_{1,2} =
\lambda_{1,2}, \lambda_{1,3} =
\frac{\lambda_{1,10} \lambda_{2,3}}{\lambda_{2,10}}, \lambda_{1,4} =
\frac{\lambda_{1,10} \lambda_{2,4}}{\lambda_{2,10}}, \lambda_{1,10} =
\lambda_{1,10}, \lambda_{2,1} = \lambda_{2,1}, \lambda_{2,2} =
\lambda_{2,2}, \lambda_{2,3} = \lambda_{2,3}, \lambda_{2,4} =
\lambda_{2,4}, \lambda_{2,10} = \lambda_{2,10}, \lambda_{3,1} =
\lambda_{3,1}, \lambda_{3,2} = \lambda_{3,2}, \lambda_{3,3} = 0,
\lambda_{3,4} = 0, \lambda_{3,10} = 0\}, \{\lambda_{1,1} =
\lambda_{1,1}, \lambda_{1,2} = \lambda_{1,2}, \lambda_{1,3} =
\lambda_{1,3}, \lambda_{1,4} = \lambda_{1,4}, \lambda_{1,10} =
\lambda_{1,10}, \lambda_{2,1} = \lambda_{2,1}, \lambda_{2,2} =
\lambda_{2,2}, \lambda_{2,3} = 0, \lambda_{2,4} = 0, \lambda_{2,10} = 0
, \lambda_{3,1} = \lambda_{3,1}, \lambda_{3,2} = \lambda_{3,2},
\lambda_{3,3} = 0, \lambda_{3,4} = 0, \lambda_{3,10} = 0\}, \{
\lambda_{1,1} = \lambda_{1,1}, \lambda_{1,2} = \lambda_{1,2},
\lambda_{1,3} = \frac{\lambda_{1,10} \lambda_{3,3}}{\lambda_{3,10}},
\lambda_{1,4} = \frac{\lambda_{1,10} \lambda_{3,4}}{\lambda_{3,10}},
\lambda_{1,10} = \lambda_{1,10}, \lambda_{2,1} = \lambda_{2,1},
\lambda_{2,2} = \lambda_{2,2}, \lambda_{2,3} =
\frac{\lambda_{2,10} \lambda_{3,3}}{\lambda_{3,10}}, \lambda_{2,4} =
\frac{\lambda_{2,10} \lambda_{3,4}}{\lambda_{3,10}}, \lambda_{2,10} =
\lambda_{2,10}, \lambda_{3,1} = \lambda_{3,1}, \lambda_{3,2} =
\lambda_{3,2}, \lambda_{3,3} = \lambda_{3,3}, \lambda_{3,4} =
\lambda_{3,4}, \lambda_{3,10} = \lambda_{3,10}\}\}$
\para

\noindent \textsf{Back to Algorithm  \ref{alg:general} at Step 9}.

\para

\noindent At this step of the algorithm one take a particular solution in the $\cV$ (see above). In first generic point of $\cV$ we take, for instance,

\vspace*{1.5mm}

\noindent $\{\lambda_{1,1} = 1, \lambda_{1,2} = -1, \lambda_{1,10} = -1,
\lambda_{2,1} = -1, \lambda_{2,2} = -3, \lambda_{2,3} = 2,
\lambda_{2,4} = -2, \lambda_{2,10} = -2, \lambda_{3,1} = -2,
\lambda_{3,2} = 0\},$

\vspace*{1.5mm}

\noindent which produces the solution
\[
\cS_a=\left(\frac{1}{2} t_{1}^{3}+\frac{1}{2} t_2 -\frac{1}{2}+\frac{1}{2} t_{2}^{3}-\frac{1}{2} t_{2}^{2}
,
t_{1}^{3}+t_{2}^{3}-t_{2}^{2}+\frac{3}{2} t_{2} +\frac{1}{2}
\right)
\]
and
\[ \begin{array}{ll}
\cQ_a=&(
64 t_{1}^{3}-96 t_{1}^{2} t_{2} +48 t_{2}^{2} t_{1} -8 t_{2}^{3}+160 t_{1}^{2}-160 t_{2} t_{1} +40 t_{2}^{2}+134 t_{1} -66 t_{2} +37
, \\&
64 t_{1}^{3}-96 t_{1}^{2} t_{2} +48 t_{2}^{2} t_{1} -8 t_{2}^{3}+160 t_{1}^{2}-160 t_{2} t_{1} +40 t_{2}^{2}+138 t_{1} -68 t_{2} +40, \\ &
-3-4 t_{1} +2 t_{2}).
\end{array}
\]
\end{example}

\section{The Computational Approach: the case of empty base locus}\label{sec:computational-approach-particual}

In this section, we present some improvements   to the previous  general computational approach. There are two sources of computational complication in Algorithm \ref{alg:general}. On the one hand, it is necessary to increase the value of the degree  $d$ until $\sol(\cP)_d\neq \emptyset$ (see Algorithm \ref{alg:subSol1})  and, on the other hand, the dimension of the linear   space can be  high.
To deal with these two difficulties we introduce two additional hypotheses. More precisely, in the sequel, we assume that the  surface $\myS$, parametrized by $\cP$, admits a birational parametrization $\cQ$ with empty base locus and we also assume that $\cP$ is transversal (see Subsection \ref{subsec:BaseLocus}). We will see that finding $\cQ$ and its corresponding $\cS$ is easier than in the general case. Indeed, Lemma \ref{lem:grado-sin-puntos-base}, provides a particular value of $d$ such that $\sol_d(\cP)\neq \emptyset$, and hence Algorithm \ref{alg:subSol1} is only executed once and, in Algorithm \ref{alg:general}, the loop in Step 3 reduces to one execution.

\para

Let us take $d$ as in Lemma \ref{lem:grado-sin-puntos-base}, and let us see that under our new hypotheses the dimension of the solution space can   be reduced.   In the general case,  $\sol_d(\cP)$ is constructed from $\myV_d$ imposing  the condition $\myF_g(\cS)=\myF_g(\cP)$.
Below we show that  $\myV_d$ can be replaced by a lower dimensional linear subsystem, generated from the base points of $\cP$ (see Subsection \ref{subsec:BaseLocus}).

\para

 \begin{theorem}\label{theoremL(D)}
 Let  $\cP$ be transversal and let $\cQ$ be a birational parametrization of $\myS$  such that $\myB(\cS)=\emptyset$. Let
 \begin{equation}\label{eq-div}
 D:=\sum_{A\in \myB(\cP)}  \sqrt{ \dfrac{\mult(A,\myB(\cP))}{\deg(\myS)}}   \, A,
 \end{equation}
and  let  $\mathscr{L}_d(D)$  be the linear system defined by the divisor $D$,  where $d$ is as in Lemma \ref{lem:grado-sin-puntos-base}. If $(\cQ,\cS)$ solves the reparametrization problem for $\cP$, where
$\cS:=(s_1:s_2:s_3)$ and $\deg(s_i)=d$, then   $s_1,s_2,s_3\in \mathscr{L}_d(D)$.
 \end{theorem}
\begin{proof} By Lemma 12, in \cite{PolinomialPerezSendra2020},    it holds that $\myB(\cS)=\myB(\cP).$
 Furthermore, by the same lemma, for $A\in \myB(\cS)$ it holds that
\begin{equation}\label{eq-multBaseS}
\mult(A,\myB(\cS))= \dfrac{\mult(A,\myB(\cP))}{\deg(\myS)}.
\end{equation}
We observe that, by Theorem 1 in \cite{PolinomialPerezSendra2020}, since $\cP$ is transversal then $\cS$ is also transversal.
Let $V_1:=x_1 s_1+x_2s_2+x_3 s_3,\, V_2:=y_1 s_1+y_2s_2+y_3 s_3$ where $x_i,y_j$ are new variables. Then, for $i\in \{1,2,3\}$ and $A\in \myB(\cS)$ it holds that
\[
\begin{array}{rclr}
\mult(A,\myC(s_i))^2 & \geq & \min\{ \mult(A,\myC(s_i))\,|\,i\in \{1,2,3\} \}^2  &  \\
& = &\mult(A,\myC(V_1))^2 & \text{(see Lemma 2 in \cite{PolinomialPerezSendra2020})}\\
&= & \mult(A,\myB(\cS)) & \text{($\cS$ is transversal)}\\
& = & \dfrac{\mult(A,\myB(\cP))}{\deg(\myS)} & \text{(see \eqref{eq-multBaseS}) }
\end{array}
\]
Therefore,
\[ \mult(A,\myC(s_i)) \geq
 \sqrt{ \dfrac{\mult(A,\myB(\cP))}{\deg(\myS)}}.  \]
Thus, since $\deg(s_i)=d$, then $s_i\in \mathscr{L}_d(D)$.
\end{proof}

 \para

 Let us define $\sol_{d}^{*}(\cP)$ as the subset of $\sol_d(\cP)$  such that if $\cS$ is induced by  $(s_1,s_2,s_3)\in\sol_d(\cP)^3$, and $(\cQ,\cS)$ solves the reparametrization problem for $\cP$, then $\myB(\cQ)=\emptyset$.
  Then the following result holds.

  \para

  \begin{corollary} With the hypotheses of Theorem \ref{theoremL(D)}, it holds that
  $\sol_{d}^{*}(\cP)\subset \mathscr{L}_d(D)$.
  \end{corollary}

  \para

 Now, one can proceed as in Subsection \ref{subsec:generalAlgoritm} but, instead of starting with $\myV_d$, we start with $\mathscr{L}_d(D)$. Therefore, the polynomial $E_d$ (see \eqref{eq-Egeneral}) is now the defining polynomial of the linear system $\mathscr{L}_d(D)$.

\para

Let us discuss how to computationally treat \eqref{eq-div}. First, let us deal with $\myB(\cP)$. We recall that (see Subsection \ref{subsec:BaseLocus} and \eqref{eq-P})
\[ \myB(\cP)=\myC(p_1)\cap \myC(p_2)\cap \myC(p_3)\cap \myC(p_4).  \]
Let us decompose $\myB(\cP)$ as
\[ \myB(\cP)=\myB(\cP)^a \cup \myB(\cP)^{\infty}, \]
where $\myB(\cP)^a$ and $\myB(\cP)^{\infty}$ represent, respectively, the sets of affine base points and base points at infinity.
Let $p(t_1,t_2):=\gcd(p_1(t_1,t_2,0), \ldots, p_4(t_1,t_2,0))$. Then
\[ \myB(\cP)^{\infty}=\{ (a:b:0)\in \mathbb{P}^{2}(\K) \,|\, p(a,b)=0\}. \]
For the affine base point, one may consider   a minimal  Gr\"obner basis $\myG_{\myB(\cP)}$, w.r.t. lexicographic ordering with $t_1<t_2$, of the ideal
$<p_1(t_1,t_2,1), \ldots,p_4(t_1,t_2,1) > \subset \K[t_1,t_2]$. Then
$\myG_{\myB(\cP)}=\{g_{1,1}, g_{2,1},\ldots,g_{2,k_2}\}
$
 where $g_{1,1}\in \K[t_1]$,
 $g_{2,j}\in \K [t_1,t_2]$  (see e.g. \cite{Winkler} page 194).  Then,
\[ \myB(\cP)^a=\{(a:b:1)\in \mathbb{P}^{2}(\K)\,|\, g(a,b)=0 \,\,\text{for $g\in \mathscr{G}_{\myB(\cP)}$}\}.\]
 Alternatively, since one has to compute the intersection of finitely many plane curves, one can approach the problem by means of resultants; see e.g. \cite{SWP}. We also observe that, after a suitable linear change of parameters $t_1,t_2$ one can proceed similarly as in Subsection \ref{subsec:generalAlgoritm} and use the Shape Lemma.

 \para

 Concerning the coefficient of the base point in the divisor in \eqref{eq-div}, we first recall that we have already commented how to determine $\deg(\myS)$.  Furthermore, the numerators, namely, $\mult(A,\myB(\cP))$ can be computed by using the curves $\myC(W_1),\myC(W_2)$ where (see \eqref{eq-W})
\[  W_1=x_1 p_1(\ot)+\cdots+ x_4 p_4(\ot),\, \,W_2=y_1 p_1(\ot)+\cdots+ y_4 p_4(\ot)
 \]
 and taking into account that
  \[ \mult(A,\myB(\cP))=\mult_A(\myC(W_1),\myC(W_2)). \]
  Moreover, since $\cP$ is transversal, then
   $$\mult(A,\myB(\cP))=\mult(A,\myC(W_1))^2.$$
 In addition, by Proposition 2.5, in \cite{CoxPerezSendra2020}, it holds that
  $$\mult(A,\myC(W_1))=\mult(A,\myC(W_2))=\min\{\mult(A,\myC(p_i))\,|\, i\in \{1,\ldots,4\}\}.$$
 Therefore, because of the transversality,
 \[ \mult(A,\myB(\cP)) =\min\{\mult(A,\myC(p_i))\,|\, i\in \{1,\ldots,4\}\}^2. \]

  In the following algorithm we use these ideas to derive a general algorithm when the additional hypotheses of this section are assumed. In order to check the transversality of the input parametrization, we refer to Algorithm 1 presented in Section 5 in \cite{PolinomialPerezSendra2020}. In addition, we assume that the algorithm in \cite{sonia} has been applied and has not provided an answer to the problem.

\para

\begin{breakablealgorithm}
\caption{\textsf{Second Main Algorithm}}
\label{alg:generalSpecial}
\begin{algorithmic}[1]
	\REQUIRE $\cP$ a transversal rational surface parametrization as in \eqref{eq-P} and assuming that there exists a birational parametrization $\cQ$ of $\myS$ such that $\myB(\cQ)=\emptyset$.
	\ENSURE    a solution $(\cQ,\cS)$ to the birational reparametrization problem for $\cP$ with $\myB(\cQ)=\emptyset$. Or a message informing that the additional hypothesis on $\cQ$ is not satisfied; in this case, one could apply Algorithm \ref{alg:general}.
	\STATE Compute $\ell={\deg(\cP)}/{\sqrt{\deg(\myS)}}$ (see \eqref{eq-degree} and Lemma \ref{lem:grado-sin-puntos-base}).
	\IF{$\ell\not\in \mathbb{N}$} \STATE \textbf{return} the additional hypothesis does not hold. \ENDIF
\STATE Compute $D$ (see \eqref{eq-div}) as well as the linear system $\mathscr{L}_d(D)$ associated to $d=\ell$.

 Let  $E_d(\Lambda,\ot)$
be the defining polynomial of $\mathscr{L}_d(D)$. Let, for $i=1,2,3$,
$E_d^i:=E(\Lambda_{i},\ot)$ where $\Lambda_i$ are different tuples of new parameters.
	\STATE  Compute $\myGG$ as in \eqref{eq-GBstar}.
	{\STATE\label{paso:subalg1}  Apply Algorithm \ref{alg:subSol1} with   $\myGG$ and $E_d^i$ computed above.
	\IF{the output of Algorithm \ref{alg:subSol1}  is the empty set}
\STATE \textbf{return} the additional hypothesis does not hold.
	\ELSE
         \STATE Let $[C,[U_1,U_2,U_3]]$ be the output of Algorithm \ref{alg:subSol1}.
         \STATE Find $\ov\in \mathbb{V}(C) \setminus \cup_{i=1}^{3} \mathbb{V}(U_i)$.
          \STATE Set $\cS:=(E_{d}^{1}(\ov):E_{d}^{2}(\ov):E_d^{3}(\ov))$.
           \FOR{$j\in \{1,2,3\}$}
         \STATE Determine $\Phi_j $ see \eqref{eq-Phi}.
         \STATE Compute the implicit equation $H_j:=A_{j,1}(x,y)-A_{j,0}(x,y)z$ of the surface parametrization $\Phi_j$.
         \ENDFOR
         \STATE Set $\cQ:=(A_{1,1}/A_{1,0}, A_{2,1}/A_{2,0},A_{3,1}/A_{3,0})$.
	 \ENDIF
	 \STATE \textbf{return} $[\cS,\cQ]$.}
	\end{algorithmic}
\end{breakablealgorithm}

\para

We finish this section with an example where we illustrate how the previous algorithm works. One may check that algorithm in \cite{sonia} does not provide an answer to the problem.

 \para

\begin{example}\label{example:general2}
Let $\cP(\ot)= (p_1: p_2: p_3: p_4)$ be the input of Algorithm \ref{alg:generalSpecial}, where\\

\noindent
$p_1(\ot)=-t_{1}^{2} \left(t_{1}^{3}-t_{1}^{2} t_{3}-2 t_{1} t_{2} t_{3}-t_{2}^{2} t_{3}\right) t_{2},$\\

\noindent
$p_2(\ot)=t_{1}^{6}-2 t_{1}^{5} t_{3}+t_{1}^{4} t_{2}^{2}-4 t_{1}^{4} t_{2} t_{3}+t_{1}^{4} t_{3}^{2}-t_{1}^{3} t_{2}^{3}+4 t_{1}^{3} t_{2} t_{3}^{2}+t_{1}^{2} t_{2}^{3} t_{3}+6 t_{1}^{2} t_{2}^{2} t_{3}^{2}+4 t_{1} t_{2}^{3} t_{3}^{2}+t_{2}^{4} t_{3}^{2}$\\

\noindent
$p_3(\ot)=t_1^4 t_2^2, $\\

\noindent
$p_4(\ot)=t_{1}^{2} \left(t_{1}^{2}+t_{1} t_{2}-t_{1} t_{3}-t_{2}^{2}\right)^{2}.$\\

$\cP$ is  a transversal   rational projective parametrization (see Subsection \ref{subsec:BaseLocus}) of an algebraic surface $\myS$ . In Step  1,  we obtain $\ell=3$. In Step 5, we compute $D$ (see \eqref{eq-div}) as well as the linear system $\mathscr{L}_d(D)$ associated to $d=\ell=3$. We get $$D=2\cdot (0:0:1)+1\cdot (1:0:1)+1\cdot (0:1:0).$$
 Let  $E_{3}^{i}(\Lambda,\ot)$
be the defining polynomial of $\mathscr{L}_D$ and for $i=1,2,3$. We get
\begin{equation}\label{eq-E3ex}
E_{3}^{i}:=E(\Lambda_{i},\ot)=t_1^3 \lambda_{i,1} + t_1^2 t_2 \lambda_{i,2} - t_1^2 t_3 \lambda_{i,1} + t_1 t_2^2 \lambda_{i,3} + t_1 t_2 t_3 \lambda_{i,4}+ t_2^2 t_3 \lambda_{i,5}
\end{equation} where $\Lambda_i$ are different tuples of new parameters. In Step 6 we compute $\myGG:=\{u_{1}^{*}(t_1), t_2v(t_1)\}$. We get
\[ \begin{array}{lcl}
u_{1}^{*}&=& \left(2 t_{1}^{2}-2 t_{1}\right) h_{1}^{4}+\left(\left(2 h_{2}-2\right) t_{1}^{2}+\left(h_{2}+2\right) t_{1}+h_{2}\right) h_{1}^{3}\\
&& +2 h_{2} \left(\left(h_{2}-\frac{5}{2}\right) t_{1}+h_{2}+\frac{1}{2}\right) t_{1} h_{1}^{2}+h_{2}^{2} t_{1}^{2} \left(h_{2}-4\right) h_{1}-h_{2}^{3} t_{1}^{2}. \\
v(t_1)&=& -\left(-2 h_{1}^{4}-2 h_{1}^{3} h_{2}-2 h_{1}^{2} h_{2}^{2}-h_{1} h_{2}^{3}+2 h_{1}^{3}+5 h_{1}^{2} h_{2}+4 h_{1} h_{2}^{2}+h_{2}^{3}\right) t_{1}^{2}\\
&&-\left(2 h_{1}^{5}+2 h_{1}^{4} h_{2}  -2 h_{1}^{4}-4 h_{1}^{3} h_{2}-2 h_{1}^{2} h_{2}^{2}\right) t_{1}+h_{1}^{4} h_{2}.
\end{array}
\]

\noindent \textsf{Algorithm  \ref{alg:subSol1}} starts.

\para

\noindent \textsf{Step 1.}
We set $C=U_1=U_2=\emptyset$.

\vspace*{1mm}

\noindent  \textsf{Step 2.}  $E_3^i(\Lambda_i, t_1, t_2),$ for $i\in \{1,2,3\}$, are as in \eqref{eq-E3ex}.

 \vspace*{1mm}

\noindent \textsf{Steps 3-6.}  $U_1=\{\lambda_{3,1},\ldots,\lambda_{3,6}\}$. For determining $U_2$, let
    \[\cE_a=\left(\frac{e_2^1}{e_2^3},\frac{e_2^2}{e_2^3}\right),\]
 where $e_2^i=E_2^i(t_1, t_2, 1)$ for $i=1,2,3$. Let   $T$ be the determinant of the Jacobian of $\cE_a$ w.r.t. $\{t_1,t_2\}$. Then, $U_2$ collects the coefficients of the numerator of $T$ w.r.t. $\{t_1,t_2\}$. In this case, $U_2$ consists in seven 3-linear polynomials in
 $\lambda_{1,i}, \lambda_{2,j}, \lambda_{3,k}$.
   \vspace*{1mm}

\noindent \textsf{Step
  7.} Let  $\tilde{G}_i^{\cE}=e_2^i(\Lambda_i, h_1, h_2)e_2^3(\Lambda_3, t_1, t_2)-e_2^3(\Lambda_3, h_1, h_2)e_2^i(\Lambda_i, t_1, t_2),\,i=1,2,3.$

   \vspace*{1mm}

\noindent \textsf{Steps 8-11.} We compute the normal form $N_1$of $\tilde{G}_1^{\cal E}$ w.r.t ${\cal G}^{\star}$, and we observe that $N_1 \not\in {\Bbb K}(h_1, h_2)[t_1, t_2]\setminus\{0\}$.

  \vspace*{1.5mm}

\noindent \textsf{Steps 12-13.} We determine $C_{1,3}$ that is
  is the set non-zero coefficients w.r.t. $\{h_1, h_2\}$ of the coefficients w.r.t.
$\{t_1, t_2\}$  of $N_1$.  We solve the system of bilinear forms $C_{1,3}$. We get $\cV$ as the set solutions:

 \vspace*{1.5mm}

  \noindent
  $\{\{\lambda_{1,1} = 0, \lambda_{1,2} = 0, \lambda_{1,3} = 0,
\lambda_{1,4} = \lambda_{1,4}, \lambda_{1,5} = 0, \lambda_{3,1} = 0,
\lambda_{3,2} = 0, \lambda_{3,3} = 0, \lambda_{3,4} = \lambda_{3,4},
\lambda_{3,5} = 0\}, \{\lambda_{1,1} = 0, \lambda_{1,2} =
\lambda_{1,2}, \lambda_{1,3} = 0, \lambda_{1,4} = \lambda_{1,4},
\lambda_{1,5} = 0, \lambda_{3,1} = 0, \lambda_{3,2} = \lambda_{3,2},
\lambda_{3,3} = 0, \lambda_{3,4} =
\frac{\lambda_{1,4} \lambda_{3,2}}{\lambda_{1,2}}, \lambda_{3,5} = 0
\}, \{\lambda_{1,1} = \lambda_{1,1}, \lambda_{1,2} =
\lambda_{1,2}, \lambda_{1,3} = 0, \lambda_{1,4} = \lambda_{1,4},
\lambda_{1,5} = 0, \lambda_{3,1} = \lambda_{3,1}, \lambda_{3,2} =
\frac{\lambda_{1,2} \lambda_{3,1}}{\lambda_{1,1}}, \lambda_{3,3} = 0,
\lambda_{3,4} = \frac{\lambda_{1,4} \lambda_{3,1}}{\lambda_{1,1}},
\lambda_{3,5} = 0\}, \{\lambda_{1,1} = \lambda_{1,1},
\lambda_{1,2} = \lambda_{1,2}, \lambda_{1,3} = \lambda_{1,3},
\lambda_{1,4} = \lambda_{1,4}, \lambda_{1,5} = 0, \lambda_{3,1} =
\frac{\lambda_{1,1} \lambda_{3,3}}{\lambda_{1,3}}, \lambda_{3,2} =
\frac{\lambda_{1,2} \lambda_{3,3}}{\lambda_{1,3}}, \lambda_{3,3} =
\lambda_{3,3}, \lambda_{3,4} =
\frac{\lambda_{1,4} \lambda_{3,3}}{\lambda_{1,3}}, \lambda_{3,5} = 0
\}, \{\lambda_{1,1} = \lambda_{1,1}, \lambda_{1,2} =
\lambda_{1,2}, \lambda_{1,3} = \lambda_{1,3}, \lambda_{1,4} =
\lambda_{1,4}, \lambda_{1,5} = \lambda_{1,5}, \lambda_{3,1} = 0,
\lambda_{3,2} = 0, \lambda_{3,3} = 0, \lambda_{3,4} = 0, \lambda_{3,5}
 = 0\}, \{\lambda_{1,1} = \lambda_{1,1}, \lambda_{1,2} = \lambda_{1,2}
, \lambda_{1,3} = -\lambda_{1,1}-\lambda_{1,5}, \lambda_{1,4} =
2 \lambda_{1,5}, \lambda_{1,5} = \lambda_{1,5}, \lambda_{3,1} =
\lambda_{3,1}, \lambda_{3,2} = \lambda_{3,2}, \lambda_{3,3} =
-\lambda_{3,1}-\lambda_{3,5}, \lambda_{3,4} = 2 \lambda_{3,5},
\lambda_{3,5} = \lambda_{3,5}\}, \{\lambda_{1,1} =
\frac{\lambda_{1,5} \lambda_{3,1}}{\lambda_{3,5}}, \lambda_{1,2} =
\frac{\lambda_{1,5} \lambda_{3,2}}{\lambda_{3,5}}, \lambda_{1,3} =
\frac{\lambda_{1,5} \lambda_{3,3}}{\lambda_{3,5}}, \lambda_{1,4} =
\frac{\lambda_{1,5} \lambda_{3,4}}{\lambda_{3,5}}, \lambda_{1,5} =
\lambda_{1,5}, \lambda_{3,1} = \lambda_{3,1}, \lambda_{3,2} =
\lambda_{3,2}, \lambda_{3,3} = \lambda_{3,3}, \lambda_{3,4} =
\lambda_{3,4}, \lambda_{3,5} = \lambda_{3,5}\}, \{\lambda_{1,1}
 = -\lambda_{1,3}, \lambda_{1,2} = \lambda_{1,2}, \lambda_{1,3} =
\lambda_{1,3}, \lambda_{1,4} = 0, \lambda_{1,5} = 0, \lambda_{3,1} =
-\lambda_{3,3}, \lambda_{3,2} = \lambda_{3,2}, \lambda_{3,3} =
\lambda_{3,3}, \lambda_{3,4} = 0, \lambda_{3,5} = 0\}, \{\lambda_{1,1}
 = -\lambda_{1,3}-\lambda_{1,5}, \lambda_{1,2} = \lambda_{1,2},
\lambda_{1,3} = \lambda_{1,3}, \lambda_{1,4} = 2 \lambda_{1,5},
\lambda_{1,5} = \lambda_{1,5}, \lambda_{3,1} = \lambda_{3,1},
\lambda_{3,2} = \lambda_{3,2}, \lambda_{3,3} = -\lambda_{3,1},
\lambda_{3,4} = 0, \lambda_{3,5} = 0\}\}.
$\\

\noindent $\cV\subset \myV_d\times \myV_d$ contains 9 generic solutions.

  \vspace*{1mm}

\noindent \textsf{Step 14-18}   We construct $\cV_{1,2,3}$. We delete from  $\cV$   those solutions  vanishing all polynomials in, at least, one of the sets $U_i,\, i =
1,2$. We get that  $\cV_{1,2,3}=\cV$.

  \vspace*{1mm}

\noindent \textsf{Step 19.} We initialize $C:=\emptyset$ and $U_3:=\emptyset$.

\vspace*{1mm}

\noindent \textsf{Steps 20-40.}  Since $\#(\cV_{1,2,3})=9$ the for-loop consists   in nine iterations. We compute the specializations $\breve{G}_{1}^{\cE}, \breve{G}_{2}^{\cE}$ as well as the resultant $\breve{R}_{1}^{\cE}$. $\deg_{t_1}(\breve{R}_{1}^{\cE})=3$. After some algebraic manipulations, we that  the new $\cV$ is

   \vspace*{1mm}

  \noindent
  $\{\{\lambda_{1,1} = -\lambda_{1,3}-\lambda_{1,5}, \lambda_{1,2} =
\lambda_{1,2}, \lambda_{1,3} = \lambda_{1,3}, \lambda_{1,4} =
2 \lambda_{1,5}, \lambda_{1,5} = \lambda_{1,5}, \lambda_{2,1} =
-\lambda_{2,3}-\lambda_{2,5}, \lambda_{2,2} = \lambda_{2,2},
\lambda_{2,3} = \lambda_{2,3}, \lambda_{2,4} = 2 \lambda_{2,5},
\lambda_{2,5} = \lambda_{2,5}, \lambda_{3,1} = \lambda_{3,1},
\lambda_{3,2} = \lambda_{3,2}, \lambda_{3,3} = -\lambda_{3,1},
\lambda_{3,4} = 0, \lambda_{3,5} = 0\},\,\,\{\lambda_{1,1} = \lambda_{1,1}, \lambda_{1,2} = \lambda_{1,2},
\lambda_{1,3} = -\lambda_{1,1}-\lambda_{1,5}, \lambda_{1,4} =
2 \lambda_{1,5}, \lambda_{1,5} = \lambda_{1,5}, \lambda_{2,1} =
\lambda_{2,1}, \lambda_{2,2} = \lambda_{2,2}, \lambda_{2,3} =
-\lambda_{2,1}-\lambda_{2,5}, \lambda_{2,4} = 2 \lambda_{2,5},
\lambda_{2,5} = \lambda_{2,5}, \lambda_{3,1} = \lambda_{3,1},
\lambda_{3,2} = \lambda_{3,2}, \lambda_{3,3} =
-\lambda_{3,1}-\lambda_{3,5}, \lambda_{3,4} = 2 \lambda_{3,5},
\lambda_{3,5} = \lambda_{3,5}\}\}.
$\\

\noindent
and thus we have  2 solutions.

\para

\noindent \textsf{Back to Algorithm  \ref{alg:generalSpecial} at Step 12}.

\para

\noindent At this step of the algorithm one take a particular solution in the $\cV$. In the second generic point of $\cV$ we take, for instance,

\vspace*{1.5mm}

\noindent $\{\lambda_{1,2} = 1, \lambda_{1,3} = 1, \lambda_{1,5} = -1,
\lambda_{2,2} = 0, \lambda_{2,3} = 0, \lambda_{2,5} = 1, \lambda_{3,1}
 = 0, \lambda_{3,2} = 1\}
,$\\

\noindent which produces the solution
\[
\cS_a=\left(\frac{t_1^2 + t_1 t_2 - t_2^2 - t_1}{t_1t_2},\,-\frac{t_1^3 - t_1^2 - 2 t_1 t_2 - t_2^2}{t_1^2t_2}\right),
\]
and
\[
\cQ_a=\left(\frac{t_2}{(t_1 + t_2 - 2)^2}, -\frac{-t_2^2 + t_1 - 2}{(t_1 + t_2 - 2)^2}, \frac{1}{(t_1 + t_2 - 2)^2}\right).
\]
\end{example}

\section{Conclusions and future work}\label{S-conclusions}

In this paper, given a surface  $\myS$, rationally parametrized by $\cP(\ot)$, we present a method (see Algorithm \ref{alg:general}) that determine a birational parametrization $\cQ(\ot)$ of $\myS$ as well as a rational map $\cS:\mathbb{P}^{2}(\K)\dasharrow \mathbb{P}^{2}(\K)$ such that $\cP(\ot)=\cQ(\cS(\ot)).$

\para

In addition, we present some improvements   to this previous  general computational approach by introducing two additional hypotheses. More precisely,  we assume that the  input surface $\myS$ admits a birational parametrization $\cQ$ with empty base locus and  that $\cP$ is transversal. As a consequence of these two hypotheses, we see that the previous general method simplifies considerably.

\para

As topics for further research on the problem, we mention some potential lines of work.
\begin{enumerate}
\item In order to simplify further the running performance of the method, one may consider computing probabilistically the fibre of the input parametrization; this could be done working on a suitable open subset as described in   \cite{PS-grado} (Algorithm-2). In this case, one should compute the fibre of the input parametrization  by choosing several
random points on the surface.
\item As mentioned in the introduction, the methods presented in this paper could be used as a general strategy to approach the problem of computing birational reparametrizations of both rational surfaces of any codimension and rational varieties of any dimension. For this purpose, one would have to develop further the theory of bases points and generic fiber. For that, one may use   the results   in \cite{PSbasePoints} and study generalizations of  the results in \cite{CoxPerezSendra2020} and \cite{PS-grado} to arbitrary rational varieties.
\item In this paper, we have not paid attention to the field required to provide the birational parametrization. An interesting question is to analyze the optimality of the field extension used in the process. And specially interesting, for the case of applications, is to analyze the problem when the ground field is $\mathbb{Q}$ or a real extension of $\mathbb{Q}$ and the output is required to be expressed over $\mathbb{R}$.
\end{enumerate}
 \para

\noindent
\textbf{Acknowledgements.} Authors partially supported by the grant PID2020-113192GB-I00 (Mathematical
Visualization: Foundations, Algorithms and Applications) from the
Spanish MICINN.\\
The first two authors  belong to the Research Group ASYNACS (Ref. CCEE2011/R34).

\end{document}